 \documentstyle{amsppt}
 \document

\def\RBS{{RBS}}
 \def\AP{{\bold A_P}}
 \def\KP{{K_P}}
 \def\g{{\frak g}}
 
 \def\k{{\frak k}}
 \def\p{{\frak p}}
 \def\q{{\frak q}}

 \def\RR{{\Bbb R}}
\def\ZZ{{\Bbb Z}}
\def\QQ{{\Bbb Q}}
\def\CC{{\Bbb C}}

\def\CE{{\Cal E}}

\def\Abar{{\overline A}}
\def\CEbar{{\overline \CE}}
\def\Ebar{{\overline E}}
\def\Mba{{\overline M}}
\def\Mbar{{M^{BS}}}
\def\Xbar{{\overline X}}
 \def\nabla{{\triangledown}}

\def\reg{{\roman{reg}}}
\def\control{{\roman{ctrl}}}
\def\IC{{\Cal{IC}}}
\def\BB{{{BB}}}
\def\BS{{{BS}}}
\def\b{{{\bullet}}}
\def\RBS{{{RBS}}}

\def\sq{{$\,\,\square$}}

\hfuzz=5pt
\baselineskip=14pt

 \define\isoarrow{{\overset\sim\to\longrightarrow}}
 \def\-{{-1}}

 \NoBlackBoxes
\magnification\magstep1

\voffset -1.5cm
\hoffset 0.7cm

\centerline{\bf On the reductive Borel-Serre compactification:}
\centerline{\bf $L^p$-cohomology
of arithmetic groups (for large $p$)}
\bigskip

\centerline{Steven Zucker\footnote{Support in part by the National Science
Foundation,
through Grant DMS9820958}}
\smallskip

\centerline{Department of Mathematics, Johns Hopkins University, Baltimore,
MD 21218
USA\footnote{e-mail address: zucker\@jhu.edu}}

\bigskip\bigskip

\centerline{\bf Introduction}
\medskip

The main purpose of this article is to give the proof of the following theorem,
as well as some applications of the result.
\proclaim {Theorem 1}
Let $M$ be the quotient of a non-compact symmetric space by an
arithmetically-defined
group of isometries, and $M^\RBS$ its reductive Borel-Serre
compactification.
Then for $p$ finite and sufficiently large there is a canonical isomorphism
$$
H^\bullet_{(p)}(M)\simeq H^\bullet (M^\RBS).
$$
\endproclaim
\noindent Here, the left-hand side is the $L^p$-cohomology of $M$ with
respect to
a (locally) invariant metric.  Though it would be more natural to allow
$p=\infty$
in Theorem 1, this is not generally possible (see (3.2.2)).  On the
other hand, there is
a natural mapping $H^\bullet_{(\infty)}(M)\to H^\bullet_{(p)}(M)$ when $p<
\infty$, because $M$ has finite volume.  The
definition of
$M^\RBS$
is recalled in (1.9).

Theorem 1 can be viewed as an analogue of the so-called Zucker
conjecture (in the case of constant coefficients), where $p=2$:
\proclaim {Theorem {\rm [L],\,[SS]}} Let $M$ be the quotient of a
Hermitian symmetric space of non-compact type by an
arithmetically-defined group of isometries, i.e., a locally
symmetric variety; let $M^\BB$ its Baily-Borel Satake
compactification.  Then there is a canonical isomorphism
$$
H^\bullet_{(2)}(M)\simeq IH_{\bold m}^\bullet (M^\BB),
$$
where the right-hand side denotes the middle intersection cohomology of
$M^\BB$.
\endproclaim
\noindent However, Theorem 1 is not nearly so difficult to prove, once one
senses that
it is true; it follows
without much
ado from the methods in [Z3] (the generalization to $L^p$, $p\ne 2$, of
those of [Z1] for $L^2$).

As far as I know, the reductive Borel-Serre compactification was first
used in
[Z1,\S 4] (where it was called $Y$). This
space, a rather direct alteration of the manifold-with-corners constructed in
[BS], was
introduced there to facilitate the study of the $L^2$-cohomology of $M$.
It
also plays a central role as the natural setting for the related {\it weighted
cohomology}
of [GHM].  It
is a principal theme that $M^\RBS$ is an important space when $M$ is
an algebraic
variety over $\CC$, despite the fact that $M^\RBS$ is {\it almost never}
an
algebro-geometric, or even complex analytic, compactification of $M$.

This work had its origin in my wanting to understand [GP].  It is convenient
to formulate
the latter before continuing with the content of this article.  Let $Y$
be a Hausdorff
topological space.  For any complex vector bundle $E$ on $Y$, one has its
Chern
classes $c_k(E)\in H^{2k}(Y,\ZZ)$.
If we further assume that $Y$ is connected, compact, stratified and oriented,
then $H_{d}
(Y,\ZZ)\simeq \ZZ$, where $d$ is the dimension of $Y$; the orientation picks
out a generator $\zeta_Y$
for this homology group, known as the fundamental class of $Y$.  We
shall henceforth
assume that $d$ is even, and we write $d=2n$.  Then, if one has positive
integers
$k_i$ for $1\le i\le \ell$ such that $\sum_i k_i = n$, one can pair $c_{k_1}(E)
\cup\dots\cup
c_{k_\ell}(E)$ with $\zeta_Y$, and obtain what is called a characteristic
number,
or Chern number, of $E$.

When $Y$ is a
$C^\infty$
manifold and $E$ is a $C^\infty$ vector bundle, the Chern classes modulo
torsion
can be constructed from any connection $\nabla$ in $E$, whereupon they get
represented, via
the de Rham theorem,
by the Chern forms $c_k(E,\nabla)$ in $H^{2k}(Y)$.  (For convenience,
we will use
and understand $\CC$-coefficients here and throughout the sequel unless it
is
specified otherwise.)  If $Y$ is compact,
the Chern
numbers can be computed by integrating $c_{k_1}(E,\nabla)\wedge\,.\,.\,.\,
\wedge
c_{k_\ell}(E,\nabla)$ over $Y$.

For stratified spaces $Y$, there is a lattice of intersection (co)homology
theories,
with variable perversity $\bold p$ as parameter, as defined by Goresky and
MacPherson
[GM1].  These range from standard cohomology as minimal object, to standard
homology as
maximal, and all coincide when $Y$ is a manifold.  They can all be defined as
cohomology
with values in some constructible
sheaf whose
restriction to the regular locus $Y^\reg$ of $Y$ is just $\CC_{Y^\reg}$.
With
mappings
going in the direction of increasing perversity, we have the basic diagram
$$\matrix
H^\bullet (Y) & \rightarrow\,\, \dots\,\,\to\quad IH_{\bold p}^\bullet (Y)
\quad
\to\,\,\dots \,\,\rightarrow & H_\bullet (Y)\\
\downarrow & \uparrow & \uparrow \\
H^\bullet (Y^\reg) & \longleftarrow \qquad H_c^\bullet (Y^\reg)\qquad
\simeq &
H_\bullet (Y^\reg).\endmatrix\tag 0.1
$$
When $Y$ is compact, $\zeta_Y$ lifts to a generator of $IH^{\bold p}_d(Y)$
for all $\bold p$.

From now on, we write $M$ for $Y^\reg$, and start to view the situation in
the opposite
way, regarding $Y$ as a topological compactification of the manifold $M$.
For any vector
bundle $E$ on $M$, and bundle extension of $E$ to $\Ebar$ on $Y$, the
functoriality
of Chern classes imply that $c_\b(\Ebar)\mapsto c_\b(E)$ under the
restriction
mapping $H^\b (Y)@>\rho >> H^\b (M)$.
One might think of this as lifting the Chern class of $E$ to the cohomology
of $Y$,
but one should be aware that $\rho$ might have non-trivial kernel, so the
lift
may depend on the choice of $\Ebar$.

The case where $E = T_M$, the tangent bundle of $M$, is quite fundamental.
Finding a
vector bundle on $Y$ that extends $T_M$ is not so natural a question when $Y$
has singularities,
and one is often inclined to forget about bundles and think instead about
just
lifting the
Chern classes.  When $Y$ is a complex algebraic variety, one considers the
{\it complex}
tangent bundle $T'_M$ of $M$.  It is shown in [M] that for constructible
$\ZZ$-valued functions $F$ on
$Y$, there is a natural assignment of Chern {\it homology} classes
$c_\b(Y;F)\in
H_\b(Y)$, such that $c_\b(Y,\bold 1)$ recovers the usual Chern classes when
$Y$ is smooth.
There has been substantial interest in lifting
these classes
to the lower intersection cohomology (as in the top row of (0.1)), best to
cohomology
(the most difficult lifting problem) for the reason mentioned earlier.

Next, take for $M$ a locally symmetric variety.  For $Y$ we might consider
any
of the interesting compactifications of $M$, which include:
$M^{BB}$, the Baily-Borel Satake compactification of $M$ as an algebraic
variety [BB];
$M_\Sigma$, the smooth toroidal compactifications of Mumford (see [Mu]);
$M^{BS}$, the Borel-Serre manifold-with-corners [BS];
$M^{RBS}$, the reductive Borel-Serre compactification.  These fit into
a diagram of compactifications:
$$\matrix
M^{BS} & @>>> &  M^{RBS}  \\
 & & \downarrow \\
M_\Sigma & @>>> & M^{BB}\endmatrix\tag 0.2
$$

If one tries to compare $M^{BS}$ and $M_\Sigma$, one sees that there is a
mapping
(of compactifications of $M$) $M^{BS}\to M_\Sigma$ only in a few cases (e.g.,
$G=SU(n,1)$).
However, by a result of Goresky and Tai [GT,\,7.3], (if $\Sigma$ is
sufficiently
fine) there are continuous mappings $M_\Sigma
 @>>> M^{RBS}$
(seldom a morphism of compactifications) such that upon inserting them in
(0.2), the obvious triangle
commutes in
the homotopy category.  One thereby gets a diagram of cohomology mappings
$$\matrix
H^\b (M) & @<<< & H^\b (M^{RBS})  \\
\uparrow & \swarrow & \uparrow \\
H^\b (M_\Sigma) & @<<< & H^\b (M^{BB});\endmatrix\tag 0.3
$$
also, the fundamental classes in $H_{d}(M_\Sigma,\ZZ)$ and $H_{d}(M^
\RBS,\ZZ)$
are mapped to $\zeta_{M^{BB}}$.

Let $E$ be a (locally) homogeneous vector bundle on $M$ (an example of which
is the holomorphic tangent
bundle $T'_M$).  There always exists an equivariant connection on $E$, whose
Chern forms are
$L^\infty$ (indeed, of constant length) with respect to the natural metric on
$M$.  In [Mu],
Mumford showed that the bundle $E$ has a so-called {\it canonical extension}
to a vector
bundle $E_\Sigma$ on $M_\Sigma$, such that these Chern forms, beyond
representing
the Chern classes of $E$ in $H^\b(M)$, actually represent the Chern classes of
$E_\Sigma$
in $H^\b(M_\Sigma)$ (see our (3.2.4)). That served the useful purpose of
placing
these
classes in a ring with Poincar\'e duality, and implied Hirzebruch
proportionality for $M$.

In [GP], Goresky and Pardon lift these classes to the cohomology of $M^{BB}$,
minimal in the lattice of
interesting compactifications of $M$, so these can be pulled back
to the
other compactifications in (0.2).  (On the other hand, the bundles do not
extend to
$M^{BB}$ in any obvious way.)  They achieve this
by constructing
another connection in $E$ (see our (5.3.3)), one that has good properties near
the singular
strata of
$M^{BB}$, using features from the work of Harris and Harris-Zucker (see
[Z5,\,App.\,B]).
With this done, the Chern forms lie in the complex of controlled differential
forms
on $M^{BB}$, whose cohomology groups give $H^\b (M^{BB})$.
In the
case of the tangent bundle, the classes map to one of the MacPherson Chern
homology
classes in $H_\b (M^{BB})$, viz., $c_\b(M^\BB; \chi_M)$, where $\chi_M$ denotes
the characteristic
(i.e., indicator) function of $M\subset M^\BB$ [GP,\,15.5].

In [GT,\,9.2], an extension of $E$ to a vector bundle $E^{RBS}$ to $M^{RBS}$
is constructed;
this does not require $M$ to be Hermitian.  There, one finds the following:
\proclaim{Conjecture A {\rm [GT,\,9.5]}}
Let $M_\Sigma \to M^{RBS}$ be any of the continuous mappings constructed in
{\rm [GT]}.  Then the
canonical extension $E_\Sigma$ is isomorphic to the pullback of $E^{RBS}$.
\endproclaim
In the absence of a proof of Conjecture A, we derive the ``topological"
analogue
of Mumford's result as a consequence of Theorem 1:
\proclaim{Theorem 2}
Let $M$ be an arithmetic quotient of a symmetric space of non-compact type.
Then the Chern
forms of an equivariant connection on $M$ represent $c_\b (E^{RBS})$ in
$H^\b (M^{RBS})$.
\endproclaim
\noindent We point out that (0.3) and Conjecture A suggest that this is more
basic
in the Hermitian setting than Mumford's result.

Goresky and Pardon predict further:
\proclaim{Conjecture B {\rm [GP]}} The Chern classes of $E^{RBS}$ are the
pullback of
the classes in $H^\b (M^{BB})$ constructed in {\rm [GP]} via the quotient
mapping
$M^{RBS}\to M^{BB}$.
\endproclaim

\noindent
Our third main result is the proof of Conjecture B.
\medskip

The material of this article is organized as follows. In \S 1 we give a
canonical
construction of the bundle $E^{RBS}$ along the lines of [BS]. We next
discuss $L^p$-cohomology,
both in general in \S 2, then on arithmetic quotients of symmetric spaces
in \S 3,
achieving a proof of Theorem 1. We make a consequent observation
in (3.3) that shows how
$L^p$-cohomology
can be used to provide definitions of mappings between topological
cohomology
groups when it is unclear how to define the mappings topologically.
In \S 4, we
treat connections and the notion of Chern forms for a natural class of
vector
bundles on stratified spaces; this allows for the proof in \S 5 of both
Theorem 2
and Conjecture B.
\bigskip

\eightpoint
This article was conceived while I was spending Academic Year 1998--99 on
sabbatical at
the Institute for Advanced Study in Princeton.  I wish to thank Mark Goresky
and John Mather for helpful discussions.

\tenpoint\baselineskip=14pt
\bigskip

\centerline{\bf 1. The Borel-Serre construction for homogeneous vector bundles}
\medskip

In this section, we make a direct analogue of the Borel-Serre construction for
the total
space of a
homogeneous vector bundle on a symmetric space, and then for any neat
arithmetic
quotient $M_\Gamma$ thereof.  It defines a natural extension of the vector
bundle to
the Borel-Serre
compactification of the space.  That the bundle extends is clear, for
attaching
of a boundary-with-corners does not change homotopy type.  Our construction
retains at
the boundary much of the group-theoretic structure.  The construction is
shown
to descend to the reductive Borel-Serre compactification $M_\Gamma^\RBS$,
reproving [GT,\,9.2].
\medskip

\noindent{\bf (1.0)} {\it Convention.}  Whenever $H$ is an algebraic group
defined over
$\QQ$, we also let $H$ denote $H(\RR)$, taken with its topology as a
real
Lie group, if there is no danger of confusion.
\medskip

\noindent{\bf (1.1)} {\it Standard notions.}  Let $G$ be a semi-simple
algebraic group over $\QQ$,
and $K$ a maximal compact subgroup of $G$, and $X = G/K$.  (Note that this
implies a choice
of basepoint for $X$, namely the point $x_0$ left fixed by $K$.)

Let $\CE = G\times_K E$ be the homogeneous vector bundle on $X$ determined
by
the representation of $K$ on the vector space $E$.  The natural projection
$\pi : \CE\to X = G\times_K \{0\}$ is induced by the projection
$E\to
\{0\}$, and is $G$-equivariant.
For $\Gamma\subset G(\QQ)$ a torsion-free arithmetic subgroup, let $M_\Gamma
= \Gamma
\backslash X$.  Then $\Gamma\backslash\CE$ is the total space of a vector
bundle
$\CE_\Gamma$ over $M_\Gamma$.  (The subscript ``$\Gamma$" was
suppressed in
the Introduction.)

If $P$ is any $\QQ$-parabolic subgroup of $G$, the action of $P$ on
$X$ is transitive.
Thus, one can also describe $\CE\to X$ as $P\times_\KP E\to P\times_\KP
\{0\}$,
where $K_P = K\cap  P$.
\medskip

\noindent{\bf (1.2)} {\it Geodesic action.} Let $U_P$ denote the unipotent
radical of
$P$, and $A_P$ the lift to $P$ associated to $x_0$ of the connected
component of the maximal
$\QQ$-split torus $Z$ of $P/U_P$.  Define the {\it geodesic action of $A_P$
on $\CE$} by the formula:
$$
a\circ (p,e) = (pa,e)\tag 1.2.1
$$
whenever $p\in P$, $e\in E$ and $a\in A_P$;  this is well-defined
because for $k\in\KP$,
$$
a\circ (pk^\-,ke) = (pk^\- a,ke) = (pak^\- , ke),\tag 1.2.2
$$
as $\AP$ and $\KP$ commute.  The
geodesic
action of $A_P$ commutes with the
action of
$P$ on $\CE$, and it projects to the geodesic action of $A_P$
on $X$
as defined in [BS,\,\S 3] (in [BS], the geodesic action is expressed
in terms of $Z$, but the definitions coincide).
\smallskip

\noindent (1.2.3) {\it Remark.}  By taking $E$ to be of dimension zero, the
construction of
Borel-Serre
can be viewed as a case of ours above.  As such, there is no real need to
recall it
separately.  Conversely, a fair though incomplete picture of our construction
can
be seen by regarding $\CE$ as simply a thickened version of $X$.
\medskip

\noindent{\bf (1.3)} {\it Corners.}  The simple roots occurring in $U_P$ set
up an isomorphism
$A_P\simeq (0,\infty)^{r(P)}$, where $r(P)$ denotes the parabolic $\QQ$-rank of
$P$.
Let $\Abar_P$ be the enlargement of $A_P$ obtained by transport of structure
from
$(0,\infty)^{r(P)}\subset (0,\infty]^{r(P)}$.  Define the {\it corner
associated
to $P$}: $\CE (P) = \CE\times_{A_P}
\Abar_P$.
There is a canonical mapping $\pi (P):\CE (P)\to X(P) = X\times_{A_P}\Abar_P$.
\smallskip

\noindent (1.3.1) {\it Remark.} Though $X(P)$ is contractible, and hence
$\CE (P)$ is
trivial, (1.2.1) {\it does not}
yield a
canonical trivialization of $\CE(P)$ over $X(P)$, because of the equivalence
relation (1.2.2) determined by $K_P$.
\smallskip

Let $\infty_P$ denote the zero-dimensional $A_P$-orbit in $\Abar_P$, which
corresponds
to $(\infty,...,\infty)\in (0,\infty]^{r(P)}$.  The {\it face of $\CE (P)$
associated to $P$} is
$$
E(P) = \CE\times_{A_P} \{\infty_P\}\simeq\CE/A_P.\tag 1.3.2
$$
It maps canonically to $X/{A_P} \simeq e(P)\subset X(P)$ (from [BS,\,5.2]).
There are
{\it geodesic projections} implicit in (1.3.2), given by the rows of the
commutative diagram
$$\CD
\CE(P) & @>\widetilde\pi_P >> & E(P)\\@VVV & & @VVV\\X(P) & @>\pi_P >> &
e(P)\endCD\tag 1.3.3
$$
\medskip

\noindent{\bf (1.4)} {\it Structure of $E(P)$.}  There is a natural
$P$-action
on $E(P)$, with $A_P$ acting trivially, projecting to the action of $P$
on
$e(P)$.
We know that $e(P)$ is homogeneous under ${{}^0\! P}$ (as in [BS,\,1.1]),
isomorphic to
$P/{A_P}$,  which contains $K_P$.  We see that $E(P)$ is
isomorphic to the
homogeneous vector bundle on $e(P)$ determined by the representation of $K_P$
on $E$.
\medskip

\noindent{\bf (1.5)} {\it Compatibility}. For $Q\subset P$, there is a
canonical
embedding of $\CE (P)$ in $\CE (Q)$, given as follows.  As in [BS,\,4.3],
write
$A_Q =
A_P\times A_{Q,P}$, with $A_{Q,P}\subset A_Q$ denoting the intersection of
the
kernels
of the simple roots for $A_P$.  Then there is an embedding
$$\multline
\CE (P) = \CE\times_{A_P} \Abar_P \simeq (\CE\times_{A_{Q,P}}A_{Q,P})
\times_{A_P}
\Abar_P \\
\subset
(\CE\times_{A_{Q,P}}\Abar_{Q,P})\times_{A_P} \Abar_P\simeq\CE\times_{A_Q}
\Abar_Q
= \CE (Q).\endmultline
$$
Moreover, this projects to $X(P)\subset X(Q)$ via $\pi (Q)$.
\medskip

\noindent{\bf (1.6)} {\it Hereditary property}.  If $Q\subset P$ again, one
can view
$E(Q)$ as part of the boundary of $E(P)$, in the same way that $e(Q)$ is
part of
the boundary of $e(P)$.  This is achieved by considering the geodesic
action of
$A_{Q,P}$ on $E(P)$ ($A_P$ acts trivially), and carrying out the analogue of
(1.3).
Thus, $E(Q)\simeq E(P)/A_{Q,P} = E(P)/A_{Q}$.
\medskip

\noindent{\bf (1.7)} {\it The bundle with corners.} Using the identifications
given in
(1.5), we recall that one puts
$$
\Xbar = \bigcup_P X(P) = \bigsqcup_P e(P),\tag 1.7.1
$$
with $P$ ranging over all parabolic subgroups of $G/\QQ$, including the
improper
one ($G$ itself).  With $\Xbar$ endowed with the weak topology from the
$X(P)$'s,
this is the manifold-with-corners construction of Borel-Serre for $X$ (see
[BS,\,\S 7]).
As such, it has a tautological stratification, with the $e(P)$'s as strata.

We likewise put $\CEbar = \bigcup_P \CE (P)$, with incidences given by (1.5),
and endow it
with the weak topology.  There is an obvious projection onto $\Xbar$.  Then
$\CEbar$
is a vector
bundle over $\Xbar$ that is stratified by the homogeneous bundles $E(P)$,
given as in (1.4).
\medskip

\noindent{\bf (1.8)} {\it Quotient by arithmetic groups.}  We can see that
$G(\QQ)$ acts
as vector
bundle automorphisms on $\CEbar$ over its action as homeomorphisms of $\Xbar$
(given in
[BS,\,7.6]); also, as it is so for $\Xbar$, the action on $\CEbar$ of any
neat
arithmetic
subgroup $\Gamma$ of $G(\QQ)$ is proper and discontinuous (cf.~[BS,\,9.3]).
Then
$\CEbar_\Gamma = \Gamma\backslash\CEbar$ is a vector bundle over $M^{BS}_
\Gamma
= \Gamma\backslash\Xbar$.

Let $\Gamma_P = \Gamma\cap P$.  The action of $\Gamma_P$ (which is contained
in ${{}^0\! P}$
of (1.4)) commutes with
the geodesic
action of $A_P$.  The faces of $\CEbar_\Gamma$ are of the form
$E'(P) =
\Gamma_P\backslash E(P)$, and are vector bundles over the faces $e'(P) =
\Gamma_P
\backslash e(P)$ of $M^{BS}_\Gamma$.  By reduction
theory [BS,\,\S 9] (but see also [Z5,(1.3)]), there is a neighborhood of
${e'(P)}$
in $M_\Gamma^{BS}$ on which geodesic projection $\pi_P$ (from (1.3.3))
descends.
The same is true for $\widetilde \pi_P$ and ${E'(P)}$ (also from (1.3.3)).
\medskip

\noindent{\bf (1.9)} {\it The reductive Borel-Serre compactification.}  We
recall the
quotient
space $X^\RBS$ of $\Xbar$.  With $\Xbar$ given as in (1.7) above, one forms
the quotient
$$
X^\RBS = \bigsqcup_P X_P,\quad\left(\text{where } X_P= U_P\backslash e(P)
\right),\tag 1.9.1
$$
where $U_P$ is, as in (1.2), the unipotent radical of $P$, and endows it with
the quotient
topology from $\Xbar$.  Because $U_Q\supset
U_P$
whenever $Q\subset P$, $X^\RBS$ is a Hausdorff space (see [Z1,(4.2)]).  There
is an induced
action of $G(\QQ)$ on $X^\RBS$, for which (1.9.1) is a $G(\QQ)$-equivariant
stratification;
$G(\QQ)$ takes the stratum $X_P$ onto that of a conjugate parabolic subgroup,
with $P(\QQ)$
preserving $X_P$.  For any arithmetic group $\Gamma\subset G(\QQ)$, one has a
quotient mapping
$$
q: M^{BS}_\Gamma\to M^\RBS_\Gamma = \Gamma\backslash X^\RBS = \,\bigsqcup_P
\widehat M_P,\tag 1.9.2
$$
with $\widehat M_P = \Gamma_P\backslash X_P$.
\medskip

\noindent{\bf (1.10)} {\it Descent of $\CEbar$ to $X^\RBS$.}  Analogous to
the
description of $\Xbar$ in (1.7), we have
$$
\CEbar = \bigcup_P \CE(P) = \bigsqcup_P E(P),\tag 1.10.1
$$
and the corresponding quotient
$$
\CE^\RBS = \bigsqcup_P\left( U_P\backslash E(P)\right).\tag 1.10.2
$$
We verify that $\CE^\RBS$ is a vector bundle on $X^\RBS$.  Since $\{X(P)\}$ is
an open
cover of $\Xbar$ (see (1.7.1)), it suffices to verify this for $\CE (P)\to
X(P)$ for each $P$ separately.

Note that $U_P$ acts on $\CE (P)$ by the formula: $u\cdot (p,e,a) = (up,e,a)$,
and this
commutes with the action of $K_P\cdot A_P$.  It follows that there is a
canonical projection
$$
U_P\backslash \CE (P)\to U_P\backslash X(P).\tag 1.10.3
$$
This gives a vector bundle on $U_P\backslash X(P)$ because $U_P\cap (K_P\cdot
A_P) = \{1\}$.

Let $X^\RBS (P)$ be the image of $X(P)$ in $X^\RBS$, and $\CE^\RBS (P)$ be
the
image of
$\CE (P)$ in $\CE^\RBS$.  These differ from (1.10.3), for the $U_P$ quotient
there is too coarse
(for instance, there are no identifications on $X$ or $\CE$ in\,\, $\CE^\RBS
\to X^\RBS$).
Rather, the pullback of (1.10.3) to $X^\RBS(P)$ {\it is}\,\, $\CE^\RBS (P)$.

When $\Gamma$ is a neat arithmetic group, $\CE^\RBS_\Gamma = \Gamma\backslash
\CE^\RBS$ is a
vector bundle on $M^\RBS_\Gamma$.  This is verified in the same manner as
(1.8).
\bigskip

\centerline{\bf 2. $L^p$-cohomology}
\medskip

By now, the notion of $L^p$-cohomology, with $1\le p\le\infty$, is rather
well-established.
The case of $p=\infty$, though, is visibly different from the case of finite
$p$, and was
neglected in [Z4].  Morally, Theorem 1 is
about
$L^\infty$-cohomology, but for technical reasons we will have to settle for
$L^p$-cohomology
for large finite $p$.  It is our first goal to prove Theorem 1.
\medskip

{\bf (2.1)} {\it Preliminaries.}  Let $M$ be a $C^\infty$ Riemannian manifold.
For any
$C^\infty$ differential form $\phi$ on $M$, its length $|\phi |$
is a
non-negative continuous function on $M$.  This determines a semi-norm:
$$
|\!|\phi |\!|_p = \cases \left (\int_M |\phi (x)|^p\,dV_M(x)\right )^{\frac
1 p}
\quad \text{if } 1\le p<\infty;\\
\roman{sup}\,\{|\phi(x) |:x\in M\}\quad \text{if } p=\infty,\endcases\tag 2.1.1
$$
where $dV_M(x)$ denotes the Riemannian volume density of $M$. One says that
$\phi$
is $L^p$ if $|\!|\phi |\!|_p$ is finite.
\proclaim {(2.1.2) Definitions}
Let $w$ be a positive continuous real-valued function on the Riemannian
manifold $M$.
\smallskip

i) The {\rm [smooth] $L^p$ de Rham complex with weight $w$} is the largest
subcomplex
of the $C^\infty$ de Rham complex of $M$ consisting of forms $\phi$ such that
$w\phi$ is $L^p$, viz.
$$
A_{(p)}^\bullet (M;w)= \{\phi\in A^\bullet (M): w\phi \text{ and }
wd\phi
\text{ are } L^p\}.\tag 2.1.2.1
$$

ii) The {\rm [smooth] $L^p$-cohomology of $M$ with weight $w$} is the
cohomology
of $A_{(p)}^\bullet (M;w)$.  It is denoted $H_{(p)}^\bullet (M;w)$.
\endproclaim

\noindent We note that in the above, there is a difference with the notation
used elsewhere:
for $p\ne\infty$, $w$ might be replaced with $w^{\frac 1 p}$ in (2.1.2.1).

When $w=1$, one drops the symbol for the weight.  Note that
the complex
depends on $w$ only through rates of the growth or decay of $w$ at infinity.
When
$M$ has finite
volume, there are inclusions $A_{(p')}^\bullet (M) \hookrightarrow A_{(p)}^
\bullet (M)$
whenever $1\le p < p'\le\infty$.  The preceding
extends to metrized local systems (cf.~[Z1,\,\S 1]).  Smooth functions are
dense in
the Banach space $L^p$ for $1\le p < \infty$, but not in $L^\infty$.

We next recall the basic properties
of $L^p$-cohomology.
Let $\Mba$ be a compact Hausdorff topological space that is a
compactification
of $M$.  One
defines a presheaf on $\Mba$ by the following rule (cf.~[Z4,\,1.9]):  to any
open
subset $V$ of $\Mba$, one assigns $A_{(p)}^\bullet (V\cap M;w)$.  Because
$\Mba$ is compact (see (2.1.5,\,ii) below), the
associated
sheaf $\Cal A_{(p)}^\bullet (\Mba ;w)$ satisfies
$$
A_{(p)}^\bullet (M;w)\isoarrow \Gamma (\Mba,\Cal A_{(p)}^\bullet
(\Mba ;w)).\tag 2.1.3
$$
It follows from the definition that whenever $q:\overline M'\to \Mba$ is a
morphism of
compactifications of $M$, one has for all $p$:
$$
q_*\Cal A_{(p)}^\bullet(\overline M';w)\simeq\Cal A_{(p)}^\bullet(\Mba;w).
\tag 2.1.4
$$

\noindent (2.1.5) {\it Remarks.} i) It is easy to see that the complex
$\Cal A_{(p)}^\bullet (\Mba ;w)$ consists
of fine sheaves if and
only if
for every covering of $\Mbar$ there is a partition of unity subordinate to
that covering
consisting of functions $f$ whose differential lies in $A^1_{(\infty)}(M)$,
i.e.,
$|df|$ is a bounded function on $M$.  Thus, (2.1.4) is for $q_*$ (as written),
not for $Rq_*$ in general.
\smallskip

ii) Note that in general, the space of global sections of $\Cal A_{(p)}^
\bullet
(M;w)$,
defined in the obvious way
(or equivalently the restriction of $\Cal A_{(p)}^\bullet (\Mba ;w)$ to
$M$) is
$A_{(p),\text{loc}}^\bullet (M;w) = A^\bullet (M)$.  Without a compact
boundary,
there is no place to store the global boundedness condition.
\smallskip

The following fact makes for a convenient simplification:
\proclaim {(2.1.6) Proposition}
Let $M$ be the interior of a Riemannian manifold-with-corners $\overline M$
(i.e., the
metric is locally extendable across the boundary).  Let $\overline{\Cal A}_{
(p)}^\bullet (\Mba ;w)$
be the sub-complex of $\Cal A_{(p)}^\bullet (\Mba;w)$ consisting of forms
that are also
smooth at the boundary of $\overline M$.  Then the inclusion
$$
\overline{\Cal A}_{(p)}^\bullet (\Mba;w)\hookrightarrow {\Cal A}_{(p)}^
\bullet (\Mba;w)
$$
is a quasi-isomorphism. \sq
\endproclaim

\noindent In other words, one can calculate $H_{(p)}^\bullet (M;w)$ using only
forms
with the nicest
behavior along $\partial\Mba$.  Moreover, $\overline{\Cal A}_{(p)}^\bullet
(\Mba;w)$
admits a simpler description; for that and the proof of (2.1.6), see (2.3.7)
and (2.3.9) below.
\medskip

{\bf (2.2)} {\it The prototype.}  We compute a simple case of
$L^p$-cohomology, one that will be useful in the
sequel.
\proclaim {(2.2.1) Proposition \rm{[Z4,\,2.1]}} Let
$\RR^+$
denote
the positive real numbers,
and $t$ the
linear coordinate from $\RR$.  For $a\in\RR$, let $w_a(t) = e^{at}$.  Then
\smallskip

i) $H^0_{(p)}(\RR^+;w_a)\simeq\cases 0 \text{ if } a>0,\\
                    \CC \text{ if } a\leq 0.\endcases$
\smallskip

ii) $H^1_{(p)}(\RR^+;w_a) = 0 \text{ for all } a\ne 0$.
\endproclaim
\demo {Proof} Again, we carry this out here only for $p=\infty$.
First, (i) is obvious: it is just an issue of whether the constant functions
satisfy
the corresponding $L^\infty$ condition.  To get started on (ii), proving that
a complex is
acyclic can be accomplished by finding a cochain homotopy operator $B$
(lowering
degrees by one), such that $\phi = dB\phi + Bd\phi$.  For the cases at hand
(1-forms
on $\RR^+$), this equation reduces to $\phi = dB\phi$.

When $a<0$, one takes
$$
B(\phi)(t) = -\int_t^\infty g(x)dx\tag 2.2.2
$$
when $\phi = g(t)dt$ (placing the basepoint at $\infty$ is legitimate, as $g$
decays
exponentially).
We need to check that (2.2.2) lies in the $L^\infty$ complex.  By hypothesis,
$$
|g(t)|\leq Cw_{-a}(t)
$$
for some constant $C$.  This implies that
$$
|B(\phi)(t)|\leq \int_t^\infty |g(x)|dx\,\,\leq \,C\!\int_t^\infty \!w_{-a}(x)
dx\sim w_{-a}(t)
$$
as $t\to\infty$.
In other words, $B(\phi)(t)w_a(t) \sim 1$, which is what we wanted to show.

When $a>0$, one takes instead
$$
B(\phi)(t) = \int_1^t g(x)dx,\tag 2.2.3
$$
and shows that $|B(\phi)|(t)\sim w_{-a}(t)$, yielding the same conclusion
about
$B(\phi)$ as before. \sq
\enddemo

\noindent (2.2.4) {\it Remark.}  One can see that for $a=0$, one is talking
about
$H^1_{(p)}(\RR^+)$, which is
not even
finite-dimensional (cf.~[Z1,\,(2.40)]); $H^1_{(\infty)}(\RR^+)$ contains the
linearly independent cohomology classes
of $t^{-\nu}dt$, for all $0\le\nu\le 1$. What was essential in the proof of
(2.2.1)
was that $w_a$ and one of its anti-derivatives had equal rates of growth or
decay when
$a\ne 0$.  That is, of course, false for $a=0$.
\medskip

{\bf (2.3)} {\it Further properties of $L^p$-cohomology.} We begin with
\medskip

\proclaim {(2.3.1) Proposition {\rm (A K\"unneth formula for
$L^p$-cohomology)}}
Let $I$ be the unit interval $[0,1]$, with the usual metric.  Then for any
Riemannian
manifold $N$ and weight $w$, the inclusion $\pi^*:\Cal A^\bullet_{(p)}(N;w)
\hookrightarrow
\Cal A^\bullet_{(p)}(I\times N;\pi^*w)$ is a quasi-isomorphism; thus
$$
H^\bullet_{(p)}(I\times N;\pi^*w)\simeq H^\bullet_{(p)}(N;w).
$$
\endproclaim

\demo {Proof} The argument is fairly standard.
The formula (2.2.3) defines an operator on forms on $I$.  Because $I$ has
finite length, one has now
$$
\phi = H\phi + dB\phi + Bd\phi,\tag 2.3.2
$$
where $H$ is---well---harmonic projection: zero on 1-forms, mean value on
0-forms.
The differential forms on a
product of two
spaces decompose according to bidegree. On $I\times N$, denote the bidegree by
$(e_I,e_N)$ (thus, for a non-zero form, $e_I\in \{0,1\}$).  The
exterior derivative on
$I\times N$ can be written as $d=d_I+\sigma_Id_N$, where $\sigma_I$ is given
by $(-1)^{e_I}$.
The operators in (2.3.2) make sense for $L^p$ forms on $I\times N$, taking,
for each $q$,
forms of bidegree $(1,q)$ to forms of bidegree $(0,q)$, and
we write
them with a subscript ``$I$"; thus, we have the identity
$$
\phi = H_I\phi + d_IB_I\phi + B_Id_I\phi.\tag 2.3.3
$$
It is clear that $B_I\phi$ is $L^p$ whenever $\phi$ is.  Note that $\sigma_I$
anticommutes
with $B_I$.  We can therefore write (2.3.3) as
$$\alignat1
\phi & = H_I\phi + dB_I\phi - \sigma_Id_NB_I\phi + B_Id\phi - B_I\sigma_Id_N
\phi\tag 2.3.4\\
& = (H_I\phi +  dB_I\phi + B_Id\phi) - (\sigma_Id_NB_I\phi + B_I\sigma_Id_N
\phi).\endalignat
$$
Since $\sigma_I$ and $d_N$ commute, the subtracted term equals $(\sigma_IB_I +
B_I\sigma_I)
d_N\phi = 0$, so (2.3.4) is just $\phi = (H_I\phi +  dB_I\phi + B_Id\phi)$.
This
implies first that $dB_I\phi$ is $L^p$ and then our assertion. \sq
\enddemo

We next use a standard smoothing argument in a neighborhood of $0\in\RR$.  To
avoid unintended
pathology, {\it we consider only monotonic weight functions} $w$.  Given
a smooth function
$\psi$ on $\RR$ of compact support, let
$$
(\Psi f)(t) = (\psi\ast f)(t) = \int\psi (x)f(t-x)dx = \int\psi
(t-x)f(x)
dx,\tag 2.3.5
$$
defined for those $t$ for which the integral makes sense.
The discussion separates into two cases:
\smallskip

i) $w(t)$ is a bounded non-decreasing function of $t$.  In this case, take
$\psi$ to be supported in $\RR^-$.
\smallskip

ii)  Likewise, when $w(t)$ blows up as $t\to 0^+$ take $\psi$ to be supported
in $\RR^+$,
and set $f(x)=0$ for $x\leq 0$.

\proclaim{(2.3.6) Lemma}
If $f\in L^p(\RR^+,w)$ (and $\psi$ is chosen as above), then $\Psi f$
is
also in $L^p(\RR^+,w)$.
\endproclaim
\demo {Proof} For $p<\infty$, see [Z4,\,1.5].  When $p=\infty$, we consider
each
of the above cases.  In case (i), we have:
$$
w(t)\Psi f(t) = \int\psi(t-x)w(t)f(x)dx = \int\psi(t-x)w(x)f(x)\{w(t)
w(x)^\-\}dx.
$$
By hypothesis, the integral involves only those $x$ for which $t<x$, and
there
$w(t)w(x)^\-\leq 1$.  It follows that $w(t)\Psi f(t)$ is uniformly
bounded.
In case (ii), when $w(t)$ blows up as $t\to 0^+$ the argument is
similar and is left to the reader. \sq
\enddemo

We use (2.3.6) to prove:
\proclaim {(2.3.7) Proposition}
With $w$ restricted as above, let $\overline{\Cal A}^\bullet_{(p)} (I;w)$
denote the
subcomplex of $\Cal A^\bullet_{(p)} (I;w)$ consisting of forms that are
smooth at
$0$.  Then the inclusion
$$
\overline{\Cal A}^\bullet_{(p)} (I;w)\hookrightarrow \Cal A^\bullet_{(p)}
(I;w)
$$
is a quasi-isomorphism, with $\Psi$ providing a homotopy inverse.
\endproclaim
\demo{Proof}
There is a well-known homotopy smoothing formula, which is at bottom a
variant
of (2.3.2).
We use the version given in [Z4,1.5], valid on the level of germs at $0$:
$$
1-\Psi = dE + Ed,\quad E = (1-\Psi )B,
$$
with $B$ as above.  Our assertions follow immediately. \sq
\enddemo

The behavior of $w$ forces the value $f(0)$ of a
function
$f\in L^\infty(I;w)\cap \overline A(I)$ to be 0 precisely in case (ii) above.
Thus we have:
\proclaim {(2.3.8) Corollary} Write $\overline I$ for the closed interval
$[0,1]$.  For the two cases preceding {\rm (2.3.6)},
$$
A^\bullet_{(\infty )} (I;w)\approx \cases A^\bullet (\overline I)\text{ \rm in
case (i),}\\
A^\bullet (\overline I,0)\text{ \rm in case (ii)}.\endcases
$$
\endproclaim

There are several standard consequences and variants of (2.3.7) in higher
dimension.
The simplest to state are (2.1.6) and its corollary; we now give the latter:

\proclaim{(2.3.9) Proposition}  Let $M$ be the interior of a Riemannian
manifold-with-corners
$\overline M$, and let $A_{(p)}^\bullet (\overline M)$ be the subcomplex of
$A_{(p)}^
\bullet (M)$
consisting of forms that are smooth at the boundary.  Then the inclusion
$$
A_{(p)}^\bullet (\overline M)\hookrightarrow A_{(p)}^\bullet (M)
$$
induces an isomorphism on cohomology.  Thus the $L^p$-cohomology of $M$
can be
computed as the cohomology of $A_{(p)}^\bullet (\overline M)$, i.e.,
$H^\bullet_{(p)}
(M)\simeq H^\bullet (M)$. \sq
\endproclaim

Finally, we will soon need the following generalization of (2.3.1):
\proclaim
{(2.3.10) Proposition}  Let $w_M$ and $w_N$ be positive
functions
on the Riemannian manifolds $M$ and $N$ respectively.  Suppose that on the
Riemannian
product $M\times N$, one has in the sense of operators on $L^p$ that $d=d_M
\otimes 1_N +
\sigma_M\otimes d_N$,
and that $H_{(p)}^\b (N;w_N)$ is finite-dimensional.  Then
$$
H_{(p)}^\b (M\times N;w_M\times w_N)\simeq H_{(p)}^\b (M;w_M)\otimes
H_{(p)}^\b
(N;w_N).
$$
\endproclaim

\noindent {\it Remarks.}  i) The condition on $M\times N$ is asserting that
the forms
on $M\times N$ that have {\it separate} $L^p$ exterior derivatives along $M$
and along $N$
are dense in the graph norm (cf.~[Z1,\,pp.178--181] for some discussion of
when
this condition holds.)
\smallskip

ii) When $p=2$, the above proposition recovers only a
special
case of what is in [Z1,\,pp.180--181]; however, the full statement of the
latter does
generalize to all values of $p$, by a parallel argument.

\demo {Proof of {\rm (2.3.10)}}  The argument is similar to what one finds in
[Z1,\S 2], which
is for the case
$p=2$, though we cannot use orthogonal projection here.  Let $h^\b = h_p^\b
(N;w_N)$
be any space of cohomology representatives for $H_{(p)}^\b
(N;w_N)$;
by hypothesis, $h^\b$ is a finite-dimensional Banach space.  It
suffices to
show that the inclusion
$$
A_{(p)}^\b (M;w_M)\otimes h_p^\b (N;w_N) \overset{\iota}\to\hookrightarrow
A_{(p)}^\b
(M\times N;w_M\times w_N)\tag 2.3.10.1
$$
induces an isomorphism on cohomology.

For each $i$, let $Z^i$ denote the closed forms in  $A^i=A_{(p)}^i (N;w_N)$.
Then $D^i=
dA^{i-1}$ is a complement to $h^i$ in $Z^i$; it is automatically closed
because of
the finite-dimensionality of $h^i$.  By the Hausdorff maximal principle, there
is a closed
linear complement $C^i$ to $Z^i$ in $A^i$ (canonical complements exist when
$p=2$).
Then the open mapping theorem of functional analysis (applied for the $L^p$
graph norm
on $A^i$), gives that the direct sum of Banach spaces,
$$
h^i\oplus D^i\oplus C^i,\tag 2.3.10.2
$$
is boundedly isomorphic to $A^i$.  With respect to this
decomposition
of $A^i$, $d_N$ breaks into the $0$-mapping on $Z^i$ and an isomorphism
$d^i: C^i\to D^{i+1}$.

We can now obtain a cochain homotopy for $A^\b$.  Let $B^i$ denote the
inverse
of $d^{i-1}$,
and $B$ and $d$ the respective direct sums of these.  One calculates that
$dB + Bd$ is
equal to $1-q$, where $q$ denotes projection onto $h^\b$ with respect to
(2.3.10.2).

Adapting this formula to $M\times N$ runs a standard course.  First, $B$
defines
an operator $B_N = 1_M\otimes B$ on $M\times N$, and likewise does $q$.  We
have the identity $1-q_N = d_N B_N + B_Nd_N$.
Noting that
$d_N$ commutes with $\sigma_M$ and that $\sigma_M^2 = 1_M$, we obtain
$$
(1-q_N) = \sigma_Md_N(\sigma_MB_N) + (\sigma_MB_N)\sigma_Md_N,
$$
and likewise $d_M(\sigma_MB_N) +(\sigma_MB_N)d_M = 0.$  Adding, we get
$1-q_N =
d\widetilde B + \widetilde B d$, with $\widetilde B = \sigma_MB_N$, and this
gives what
we wanted to know about (2.3.10.1), so we are done. \sq
\enddemo
\bigskip

\centerline{\bf 3. $L^p$ cohomology on the reductive Borel-Serre
compactification}
\medskip

In this section, we determine the cohomology sheaves of $\Cal A_{(p)}^\bullet
(M^\RBS)$
for large finite values of $p$, and compare the outcome to that of related
calculations.
\medskip

{\bf (3.1)} {\it Calculations for $M^\RBS$, and the proof of Theorem 1.} We
first observe
that $\Cal A_{(p)}^\bullet (M^\RBS)$ is a complex of fine sheaves, for the
criterion of
(2.1.5,\,i) was verified in [Z1].  (The analogous statement on $\Mbar$
is false
unless $M$ is already compact; indeed, this is why the space $M^\RBS$ was
introduced.)

Let $y\in U_P\backslash e(P)\subset M^\RBS$.  The issue is local in
nature, so
it suffices to work with $\widetilde q: M^{BS}_{\Gamma_{U_P}} \to X^\RBS,$
and therefore
we lift $y$ to $\widetilde y\in X^\RBS$.  The fiber $\widetilde q^\-
(\widetilde y)$
is the compact nilmanifold $N_P = \Gamma_{U_P}\backslash U_P$.  Since $N_P$
is compact,
neighborhoods of $\widetilde y$ in $X^\RBS$ give, via $\widetilde q^\-$, a
fundamental
system of neighborhoods of $N_P$ in $M^{BS}_{\Gamma_{U_P}}$.

As in [Z1,(3.6)], the intersection with $M_{\Gamma_{U_P}}$ of such a
neighborhood
is of the form
$$
A_P^+\times V \times N_P,\tag 3.1.1
$$
where $A_P^+\simeq (\RR^+)^{r(P)}$
and $V$
is a coordinate cell on $\widehat M_P$ (notation as in (1.9)).  After taking the
exponential
of the $A_P^+$-variable, the metric is given, up to quasi-isometry, as
$$
\sum_i dt_i^2 + dv^2 + \sum_\alpha e^{-2\alpha}du_\alpha,\tag 3.1.2
$$
where $\alpha$ runs over the roots in $U_P$.  By the K\"unneth formula
(2.3.1),
we may replace $V$
by a point in (3.1.2); we are reduced to determining $H^\bullet_{(p)}(A_P^+
\times N_P)$,
where the metric is $\sum_i dt_i^2 + \sum_\alpha e^{-2\alpha}du_\alpha$.

The means of computing this runs parallel to the discussion in [Z1,(4.20)].
We consider
the inclusions of complexes
$$\multline
\bigoplus_\beta \left( A^\bullet_{(p)}(A_P^+;w_\beta)\otimes H^\bullet_\beta
(\frak u_P,\CC)\right)
\hookrightarrow
\bigoplus_\beta \left( A^\bullet_{(p)}(A_P^+;w_\beta)\otimes \wedge^\bullet_
\beta
(\frak u_P)^*\right)\\
\simeq A^\bullet_{(p)}(A_P^+\times N_P)^{U_P}\hookrightarrow A^
\bullet_{(p)}(A_P^+\times N_P).\endmultline\tag 3.1.3
$$
Here, $\frak u_P$ denotes the Lie algebra of $U_P$, and
$$
w_\beta (a) = a^{p\beta}a^{-\delta} = a^{p\beta -\delta} = a^{p(\beta -\frac
\delta p)}
\qquad (a_i = e^{t_i}),\tag 3.1.4
$$
where $\delta$ denotes the sum of the positive $\QQ$-roots
(cf.~(3.1.9,\,ii) below).
We can see that the contribution of $\delta$ (which enters because of the
weighting of the
volume form of $N_P$) is non-zero yet increasingly negligible as $p\to\infty$.

The second inclusion in (3.1.3) is that of the ``$U_P$-invariant" forms.
Note that
this reduces considerations on $N_P$ to a finite-dimensional vector space,
viz.~$\wedge^
\bullet (\frak u_P)^*$.  Here, one is invoking the isomorphism
$$
H^\bullet (N_P)\simeq H^\bullet (\frak u_P,\,\CC)\tag 3.1.5
$$
for nilmanifolds, which is a theorem of Nomizu [N].
The exterior algebra
decomposes into non-positive weight spaces for $\frak a_P$, which we write as
$$
\wedge^\bullet (\frak u_P)^* = \bigoplus_\beta \wedge^\bullet_\beta
(\frak u_P)^*.
\tag 3.1.6
$$
The first inclusion in (3.1.3) is given by Kostant's embedding [K,(5.7.4)]
of $H^\bullet
(\frak u_P,\CC)$ in $\wedge^\bullet (\frak u_P)^*$ as a set of cohomology
representatives,
and it respects $\frak a_P$ weights.  Our main $L^p$-cohomology computation is
based on:

\proclaim {(3.1.7) Proposition} For all $p\ge 1$, the inclusions in {\rm
(3.1.3)}
are quasi-isomorphisms.
\endproclaim
\demo {Proof}
This is asserted in [Z1,(4.23),(4.25)] for the case $p=2$.  The proof given
there was
presented with $p=2$ in mind, though there is no special role of $L^2$ in it
(cf.~[Z3,(8.6)]).
We point out that [Z1,(4.25)] is about the finite-dimensional linear algebra
described
above, and that the proof (4.23) of [Z1] goes through because the process of
averaging
a function over a circle (hence a nilmanifold, by iteration) is bounded in
$L^p$-norm.
As such, one sees rather easily that the proof carries over verbatim for
general
$p$, and (3.1.7) is thereby proved.  \sq
\enddemo

\demo{Remark} We wish to point out and rectify a small mistake in the
argument
in [Z1,\,\S 4], one that ``corrects itself".  It is asserted that the
second
terms in (4.37) and (4.41) there vanish {\it by $U_{j-1}$-invariance}.  This
is false in
general.  However, the two expressions actually differ only by a sign, and
they cancel,
yielding the conclusion of (4.41).
\enddemo

We next show how (3.1.7) yields the determination of
$H^\bullet_{(p)}
(A_P^+\times N_P)$.  We may use the first complex in (3.1.3) for this purpose.
The weights in
(3.1.6) are non-positive, and (3.1.4) shows that once $p$ is sufficiently
large,
$w_\beta$ blows up exponentially in some direction whenever $\beta\ne 0$, and
decays exponentially
when $\beta = 0$.  Applying (2.2.1) and the K\"unneth theorem, we obtain:

\proclaim {(3.1.8) Corollary} For sufficiently large $p<\infty$,
$H^\bullet_{(p)}
(A_P^+\times N_P)
\simeq H^0(\frak u_P,\CC)\simeq \CC$.
\endproclaim

\noindent (3.1.9) {\it Remark}.  i) We can specify what ``sufficiently large''
means, using
(3.1.4).  Write $\delta$ as a (non-negative) linear combination of the simple
$\QQ$-roots:
$\delta = \sum_\beta\,c_\beta\beta$.  Then we mean to take $p > \max\{c_\beta
\}$.
\smallskip

ii) When $p=\infty$, one runs into trouble with
the
infinite-dimensionality of the unweighted $H^1_{(\infty)}(\RR^+)$ (see
(2.2.4)).
By using instead large finite $p$, we effect a perturbation away from the
trivial weight,
thereby circumventing the problem.
\smallskip

There is a straightforward globalization of (3.1.8), which we now state:
\proclaim {(3.1.10) Theorem}
For sufficiently large $p$, the inclusion
$$
\CC_{M^\RBS}\to q_*\Cal A_{(p)}^\bullet (\Mbar)\simeq\Cal
A_{(p)}
^\bullet (M^\RBS)
$$
is a quasi-isomorphism. \sq
\endproclaim
From this follows Theorem 1:
\proclaim {(3.1.11) Corollary} For sufficiently large finite $p$, $H_{(p)}^
\bullet(M)
\simeq H^\bullet (M^\RBS )$.
\endproclaim

{\bf (3.2)} {\it An example (with enhancement).}  Take first $G = SL(2)$.
Then $M$
is a modular curve.  There are only two distinct interesting
compactifications
(those in
(0.2)): one is $M^\BS$, and the other is $M^\RBS$ (which is homeomorphic to
$M^\BB$ and
$M_\Sigma$).  A deleted neighborhood of a boundary point (cusp) of $M^\BB$
is a
Poincar\'e punctured disc $\Delta^*_R = \{z\in \CC : 0 < |z| < R\}$, with
$R < 1$,
with metric given in polar coordinates by $ds^2 = (r\,|\log r|)^{-2}(dr^2 +
(rd\theta)^2)$.
Because $R < 1$, the metric is smooth along the boundary circle $|z|=R$.
Setting
$u = \log |\log r|$ converts the metric to $ds^2 = du^2 + e^{-2u}d\theta^2$
(recall
(3.1.3)).  One obtains from (2.3.1) and (3.1.7):
\proclaim
{(3.2.1) Proposition} Write $\Delta_R$ for $\Delta^*_R\cup \{0\}$.  Then for
the
Poincar\'e metric on $\Delta^*_R$,
$$
H^\b_{(p)} (\Delta^*_R)\simeq H^\b (\Delta_R)\simeq\CC\quad \text{ whenever }
1< p < \infty.\quad\square
$$
\endproclaim
\noindent (3.2.2) {\it Remark.}  When $p =1$, $H^2_{(1)}(\Delta^*_R)$ is
infinite-dimensional,
as is $H^1_{(\infty)}(\Delta^*_R)$; this follows from (2.2.4) and (3.1.7).
By using a
Mayer-Vietoris
argument, in the same manner as [Z1,\S 5], we get that when $M$ is a modular
curve, we
see that $H^1_{(\infty)}(M)$ is likewise infinite-dimensional.  Thus the
assertion in
(3.1.11) fails to hold for $p=\infty$, already when $G=SL(2)$.
\medskip

Using the K\"unneth formula (2.3.10), it is easy to obtain the corresponding
assertion
for $(\Delta^*_R)^n$:
\proclaim
{(3.2.3) Corollary} For the Poincar\'e metric on $(\Delta^*_R)^n$,
$$
H^\b_{(p)} ((\Delta^*_R)^n)\simeq H^\b (\Delta_R^n)\simeq\CC\quad \text
{ whenever }
1 < p < \infty. \quad\square
$$
\endproclaim

Now, let $M$ be an arbitrary locally symmetric variety.  The
smooth toroidal
compactifications $M_\Sigma$ are constructed so that they are complex
manifolds
and the
boundary is a divisor with normal crossings on $M_\Sigma$.  The local
pictures
of $M\hookrightarrow
M_\Sigma$ are $(\Delta^*)^k\times\Delta^{n-k}\hookrightarrow\Delta^n$, for
$0\le k\le n$.  The
invariant
metric of $M$ is usually not Poincar\'e in these coordinates, not
even asymptotically.
However, it is easy to construct other metrics which are.  We will use a
subscript
``P'' to indicate that one is using such a metric instead of the invariant
one.
We note that
such a metric depends on the choice of toroidal compactification.  The global
version of (3.2.3)
follows by standard sheaf theory:
\proclaim
{(3.2.4) Proposition} For a metric on $M$ that is Poincar\'e with respect to
$M_\Sigma$,
$$
H^\b_{(p),\roman P}(M)\simeq H^\b (M_\Sigma)\quad \text{ whenever } 1 < p <
\infty. \quad\square
$$
\endproclaim

The above proposition actually gives a reinterpretation of the method in [Mu].
There,
Mumford decided to work in the rather large complex of currents that also
gives the
cohomology of $M_\Sigma$.  However, he shows that the connection and Chern
forms involved
are ``of Poincar\'e growth", and that is equivalent to saying that they are
$L^\infty$
with respect to any metric that is asymptotically Poincar\'e near the boundary
of $M_\Sigma$.
One thereby sees that his argument for comparing Chern forms ([Mu,\,p.243],
based on
(4.3.4) below) in the complex of currents actually takes place in the
subcomplex
of Poincar\'e $L^\infty$ forms on $M$.

\demo{Remark} For
a convenient exposition of the growth estimates in the latter, see
[HZ2,(2.6)].
Since there is in general no morphism of compactifications between $M^\RBS$
and
$M_\Sigma$, the reader is warned that the
comparison
of their boundaries is a bit tricky (see [HZ1,(1.5),(2.7)] and [HZ2,(2.5)]).
\enddemo
\medskip

{\bf (3.3)} {\it On defining morphisms via $L^p$-cohomology.} We give next an
interesting
consequence of Theorem 1.  The space
$M$
has finite volume, so there is a canonical morphism (see (2.1))
$$
H^\b_{(p)}(M)\to H^\b_{(2)}(M)\tag 3.3.1
$$
whenever $p>2$.  For $p$ sufficiently large, the left-hand side of (3.3.1) is
naturally
isomorphic to $H^\b(M^\RBS)$.  For $p=2$, there is an analogous
assertion:
by the Zucker conjecture, proved in [L] and [SS] (see [Z2]), the right-hand
side is
naturally isomorphic to $IH_{\bold m}^\b (M^\BB)$, intersection
cohomology
with middle perversity $\bold m$ of [GM1].  These facts transform (3.3.1) into
the diagram
$$\matrix
H^\b(M^\RBS) &&\\
\uparrow & \searrow & \\
H^\b(M^\BB) & \rightarrow & IH_{\bold m}^\b (M^\BB).\endmatrix
\tag 3.3.2
$$
In other words,
\proclaim{(3.3.3) Proposition} The mapping in {\rm (3.3.1)} defines a
factorization
of the canonical mapping
$$
H^\b(M^\BB)\longrightarrow IH_{\bold m}^\b (M^\BB)
$$
through $H^\b(M^\RBS)$.
\endproclaim

A related assertion had been conjectured by Goresky-MacPherson and Rapoport,
and was proved
recently by Saper:
\proclaim
{(3.3.4) Proposition} Let $h:M^\RBS\to M^\BB$ be the canonical
quotient
mapping.  Then there is a quasi-isomorphism
$$
Rh_*\IC^\b_{\bold m}(M^\RBS,\QQ)\approx \IC^\b_{\bold m}
(M^\BB,\QQ),
$$
where $\IC$ denotes sheaves of intersection cochains.
\endproclaim

This globalizes to an isomorphism $IH^\b_{\bold m}(M^\RBS,\QQ)\isoarrow
IH^\b_{
\bold m}(M^\BB,\QQ)$, which underlies (3.3.3), enlarging the triangle
into a
commutative square defined over $\QQ$:
$$\CD
H^\b(M^\RBS)&@>>> &IH_{\bold m}^\b(M^\RBS) \\
@AAA & &@VV\simeq V \\
H^\b(M^\BB) & @>>> & IH_{\bold m}^\b (M^\BB).\endCD
$$
\newpage

\centerline{\bf 4. Chern forms for vector bundles on stratified spaces}
\medskip

In this section, we will treat the de Rham theory for stratified
spaces
that will be needed for the proof of Theorem 2.  We also develop the
associated
treatment of Chern classes for vector bundles.
\medskip

{\bf (4.1)} {\it Differential forms on stratified spaces.}  Let $Y$ be a
paracompact space
with an abstract
prestratification (in the sense of Mather) by $C^\infty$ manifolds.
Let $\Cal S$
denote the set of strata of $Y$.  If $S$ and $T$ are strata, one writes
$T\prec S$
whenever $T\ne S$ and $T$ lies in the closure $\overline S$ of $S$.

The notion of a prestratification specifies a system $\Cal C$ of {\it
Thom-Mather
control data} (see [GM2,\,p.\,42],[V1],[V2]), and
that entails
the following.  For each stratum $S$ of $Y$, there is a neighborhood $N_S$
of
$S$ in $Y$, a retraction $\pi_S:N_S\to S$,
and a continuous ``distance function" $\rho_S:N_S\to [0,\infty)$ such that
$\rho_S^\- (0)
= S$, subject to:
\medskip
{\bf (4.1.1) Conditions.}
{\it Whenever $T\preceq S$, put $N_{T,S} = N_T\cap S$, $\pi_{T,S} = \pi_T|_{
N_{
T,S}}$, and $\rho_{T,S} = \rho_T|_{N_{T,S}}$.  Then:
\smallskip

i) $\pi_T(y) = \pi_{T,S}(\pi_S(y))$ whenever both sides are defined, viz., for
$y\in N_T
\cap(\pi_S)^\-N_{T,S}$; likewise $\rho_T(y) = \rho_{T,S}(\pi_S(y))$.
\smallskip

ii) The restricted mapping $\pi_{T,S}\times\rho_{T,S}:N^\circ_{T,S}\to T
\times\RR^+$,
where $N^\circ_{T,S}=  N_{T,S}-T$,
is a $C^\infty\!$ submersion.}
\medskip
\noindent (The above conditions will be relaxed after (4.1.4) below.) A
prestratified
space is, thus, the triple $(Y,\Cal S,\Cal C)$.

Let $Y^\circ$ denote the open stratum of $Y$.  One understands that
when $S=Y^\circ$,
one has $N_S=Y^\circ$, $\pi_S = \bold 1_{Y^\circ}$ and $\rho_S\equiv 0$.
From (4.1.1), it follows
that
for all $S\in\Cal S$, $\pi_{T,S}|_{N^\circ_{T,S}}$ is a submersion; moreover,
the
closure $\overline S$ of $S$ in $Y$ is stratified by $\{T\in \Cal S: T
\preceq S\}$,
and $\Cal C_S = \{(\pi_{T,S},\rho_{T,S}): T\prec S\}$ is a system of control
data for $\overline S$.

We also recall the following (see [V2,\,Def.\,1.4]):
\proclaim{(4.1.2) Definition} A {\rm controlled mapping} of prestratified
spaces,
$\roman f:(Y,\Cal S,\Cal C)\to (Y',\Cal S',\Cal C')$, is a continuous mapping
$f:Y\to Y'$ satisfying:
\smallskip

i) If $S\in \Cal S$, there is $S'\in\Cal S'$ such that $f(S)\subseteq S'$, and
moreover,
$f|_S$ is a smooth mapping of manifolds.
\smallskip

ii) For $S$ and $S'$ as above, $f\circ \pi_S = \pi_{S'}\circ f$ in a
neighborhood of $S$.
\smallskip

\smallskip

iii) For $S$ and $S'$ as above, $\rho_{S'}\circ f = \rho_S$ in a
neighborhood of $S$.
\endproclaim

Let $j:Y^\circ\hookrightarrow Y$ denote the inclusion.
A subsheaf $\Cal A^\bullet_{Y,\Cal C}$
of
$j_*A^\bullet_{Y^\circ}$, the complex of {\it $\Cal C$-controlled $\CC$-valued
differential
forms} on $Y$, is the sheafification of the following presheaf: for $V$ open
in $Y$, put
$$
A^\bullet_{Y,\Cal C}(V) = \{ \varphi\in A^\bullet (V\cap Y^\circ): \varphi|_{
N_S\cap V\cap
Y^\circ}\in \text{ im}\,\pi_S^*\text{ for all }S\in\Cal S\}.\tag 4.1.3
$$

\noindent (4.1.4) {\it Remark.} From (4.1.1,\,i), one concludes that the
condition
in (4.1.3) for $T$ implies the same for $S$ whenever $S\succ T$, as $(\pi_T)^*
\varphi =
(\pi_S)^*(\pi_{T,S})^*\varphi$.
\smallskip

We observe that the definition of $\Cal A^\bullet_{Y,\Cal C}$ is
independent
of the distance functions $\rho_S$.  Indeed, all that we will
need from the
control data for most purposes is the collection of germs of $\pi_S$ along $S$.
We term
this {\it weak control data} (these are the equivalence classes implicit in
[V2,\,Def.\,1.3]).
In this spirit, one has the notion of a {\it weakly controlled mapping},
obtained from
(4.1.2) by discarding item (iii); cf.~(5.2.2).
The main role that $\rho_S$ plays here is to specify a model for the {\it
link} of $S$:
$$
L_S = \pi_S^\-(s_0)\cap\rho_S^\- (\varepsilon )\tag 4.1.5
$$
for any $s_0\in V_S$ and sufficiently small $\varepsilon > 0$, but the link is
also independent
of $\Cal C$; besides, we will not need that notion in this paper.
\medskip

The following is well-known:
\proclaim
{(4.1.6) Lemma} Let $Y$ be a space with prestratification.  For any open
covering
$\frak V$ of $Y$, there is a partition of unity $\{ f_V: V\in\frak V\}$
subordinate to
$\frak V$ that consists of $\Cal C$-controlled functions. \sq
\endproclaim
\noindent This is used in [V1,\,p.\,887] to prove the stratified version of
the de Rham theorem:
\proclaim
{(4.1.7) Proposition} Let $\Cal C$ be a system of (weak) control data on $Y$.
Then the
complex $\Cal A^\bullet_{Y,\Cal C}$ is a fine resolution of the constant
sheaf $\,\CC_Y\!$. \sq
\endproclaim
\proclaim
{(4.1.8) Corollary} A closed $\Cal C$-controlled differential form on $Y$
determines
an element of $H^\b(Y)$. \sq
\endproclaim
\medskip

{\bf (4.2)} {\it Controlled vector bundles.}
We start with a basic notion.
\proclaim {(4.2.1) Definition}  A {\rm $\Cal C$-controlled vector bundle} on
$Y$ is a
topological vector bundle $E$, given with local trivializations for all $V$ in
some open
covering $\frak V$ of $Y$, such that the entries of the transition matrices
are $\Cal C$-controlled.
\endproclaim

It follows from the definition that a $\Cal C$-controlled vector bundle
determines a Cech
1-cocycle
for $\frak V$ with coefficients in $GL(r,\Cal A^0_{Y,\Cal C})$.  It thereby
yields a
cohomology class in $H^1(Y,GL(r,\Cal A^0_{Y,\Cal C}))$.  The latter has a
natural interpretation:
\proclaim{(4.2.2) Proposition} The set $H^1(Y,GL(r,\Cal A^0_{Y,\Cal C}))$
is in canonical
one-to-one correspondence with the set of isomorphism classes of vector
bundles
$E$ of rank $r$ on $Y$ with $E_S=E|_S$ smooth for all $S\in\Cal S$, together
with a system
$\{\phi_S:S\in\Cal S\}$, $\{\phi_{T,S}:S,T\in\Cal S\}$ of germs of isomorphisms
of vector bundles (total
spaces) along each $T\in\Cal S$:
$$\alignat1
i)\quad & \phi_T:(\pi_T)^*E_T = E_T\times_T N_T\isoarrow E|_{N_T},
\tag 4.2.2.1\\
ii)\quad & \phi_{T,S}:
(\pi_{T,S})^*E_T = E_T\times_T N_{T,S}\isoarrow E|_{N_{T,S}}\endalignat
$$
whenever $T\prec S$, satisfying the compatibility conditions $\phi_T
= \phi_S
\circ\phi_{T,S}$.
\endproclaim
\noindent (4.2.2.2) {\it Remark}.  Condition (ii) above is, of course, the
restriction
of (i) along $S$.
\medskip

If we use $\{E_T:T\in\Cal S\}$, the stratification of $E$
induced by
$\Cal S$, then the natural projection $E_T\times_T N_T\to E_T$ gives weak
control data
for $E$.  Thus we obtain from (4.2.2):
\proclaim
{(4.2.3) Corollary}  A vector bundle $E$ on $Y$ is $\Cal C$-controlled {\rm
(as in (4.2.1))}
if and only if $E$ admits weak control data such that the bundle projection
$E\to Y$ is a
weakly controlled mapping {\rm (as in (4.1))}. \sq
\endproclaim
\demo{Proof of {\rm (4.2.2)}} Let $\xi$ be a 1-cocycle for the open covering
$\frak V$ of $Y$,
with coefficients in $GL(r,\Cal A^0_{Y,\Cal C})$. Since the functions in
$\Cal A^0_{Y,\Cal C}$
are continuous, $\xi$ determines a vector bundle of rank $r$ in the usual way;
putting $E_0$ for
$\CC^r$, one takes
$$
E=E_\xi=\bigsqcup\{(E_0\times V_\alpha): V_\alpha\in\frak V\}
$$
modulo the identifications on $V_{\alpha\beta} = V_\alpha\cap V_\beta$:
$$\matrix
E_0\times V_{\alpha\beta} & \hookrightarrow &  E_0\times V_\alpha\\
@V \bold 1\times\xi_{\alpha\beta} VV & & \\
E_0\times V_{\alpha\beta} & \hookrightarrow & E_0\times V_\beta\endmatrix\tag
4.2.2.3
$$
It is a tautology that there exist isomorphisms (4.2.2.1) locally on the
respective
bases ($Y^\circ$ or $S$), but we want it to be specified globally.

Next, let
$$
\frak V_S = \{V\in\frak V: V\cap S\ne \emptyset\},\quad \frak V(S) = \{ V\cap
S: V\in
\frak V_S\}.\tag 4.2.2.4
$$
Then $\frak V(S)$ is an open cover of $S$.  By refining $\frak V$, we may
assume without
loss of generality that $\xi_{\alpha\beta}\in\text{im}\,(\pi_S)^*$ on $V_{
\alpha\beta}
\cap N_S$ whenever $V_\alpha,V_\beta\in \frak V_S$, and write $\xi_{\alpha
\beta}=
(\pi_S)^*\,\xi^S_{\alpha\beta}$. The bundle $E_S=E|_S$ is constructed from the
1-cocycle $\xi^S$.  Let
$$
N'_S = N_S\cap\bigcup\,\{V: V\in\frak V_S\}.
$$
The relation $\xi=(\pi_S)^*\xi^S$ on $N'_S$ determines a canonical isomorphism
$\phi_S:E|_{N'_S}
\isoarrow(\pi_S)^*E_S$, for the local ones patch together; it is smooth on
each stratum
$R\succ S$.
One produces $\phi_{S,T}$ by doing the above for the restriction
of $E$ to
$\overline S$, along its stratum $T$.

The consistency condition, $\phi_T = \phi_S\circ\phi_{T,S}$ whenever $T\prec
S$, holds
because of (4.1.4).  Replacing $\frak V$ by any refinement of it, only serves
to make $N'_S$
smaller, so the germs of the pullback relations do not change.  Also, we must
check that the
isomorphisms
above remain unchanged when we replace $\xi$
by an equivalent cocycle.  Let $\xi'_{\alpha\beta}=\psi_\beta\xi_{\alpha
\beta}\psi^\-_\alpha$,
where $\psi$ is a 0-cochain for $\frak V$ with coefficients in $GL(r,\Cal
A^0_{Y,\Cal C})$.
Without loss of generality again, we assume that $\psi$ is of the form
$(\pi_S)^*\psi^S$
on $N'_S$.
The isomorphism $E(\xi'_S)\simeq E(\xi_S)$ induced by $\psi^S$ then pulls back
to the
same for the restrictions of $E(\xi')$ and $E(\xi)$ to $N'_S$, respecting the
compatibilities.

Thus, we have constructed a well-defined mapping from $H^1(Y,GL(r,\Cal A^0_{
Y,\Cal C}))$
to isomorphism classes of bundles on $Y$ with pullback data along the strata.
We wish to show that it is a bijection.

Actually, we can invert the above construction explicitly.  Given $E$,
$\phi_T$,
etc., as in (4.2.2.1), let, for each $T\in\Cal S$, $\frak V_T$ be a covering
of $T$ that gives a
1-cocycle
$\xi^T$ for $E_T$ (as a smooth vector bundle on $T$); $N'_T$ a neighborhood
of $T$, contained in $N_T$, on which the
isomorphisms
$\phi_T$ and $\phi_{T,S}$ (for all $S\succ T$) are defined; $\frak V(T) =
\pi_T^\-\frak V_T$
the corresponding covering of $N'_S$, on which $(\pi_T)^*\xi^T$ is a cocycle
giving
$E|_{N'_T}$.  Then
$$
\frak V = \bigcup\,\{\frak V_T:T\in\Cal S\}
$$
is a covering of $Y$, such that for all $V\in\frak V$, $E|_V$ has been
trivialized.

We claim that the 1-cocycle for $E$, with respect to these trivializations,
has coefficients
in $GL(r,\Cal A^0_{Y,\Cal C})$.  For $V_\alpha$ and $V_\beta$ in the same
$\frak V_T$,
we have seen already that $\xi_{\alpha\beta}$
is in $\text{im}(\pi_T)^*$.  Suppose, then, that $T\prec S$, and that
$V_\alpha\in
\frak V_T$ and $V_\beta\in\frak V_S$ have non-empty intersection.  Then
$\xi_{\alpha\beta}$
is actually in $\text{im}(\pi_S)^*$, which one sees is a consequence of the
compatibility
conditions for (4.2.2.1), and our claim is verified.

That we have described the inverse construction is easy to verify. \sq
\medskip

{\bf (4.3)} {\it Controlled connections on vector bundles.} When we speak of a
connection
on a smooth vector bundle over a manifold, and write the symbol $\nabla$ for
it, we mean
foremost the covariant derivative.  Then, the difference of two connections
is a 0-th
order operator, given by the difference of their connection matrices with
respect to any one frame.

We can define the notion of a connection on a $\Cal C$-controlled vector
bundle:
\proclaim {(4.3.1) Definition}
Let $E$ be a $\Cal C$-controlled vector bundle on $Y$.  A {\rm
$\Cal C$-controlled
connection} on $E$ is a connection $\nabla$ on $E|_{Y^\circ}$ for which there
is a covering
$\frak V$ of $Y$ such that for each $V\in\frak V$, there is a frame of
$E|_V$ such
that the connection forms lie in $\Cal A^1_{V,\Cal C}\otimes\roman{End}(E)$.
\endproclaim
\demo{Remark {\rm [added]}} It is more graceful to define a
controlled connection so as to be in accordance with (4.2.2.1): it
is a system of connections $\{(E_T,\nabla_T)\}$, with germs of
isomorphisms
$$
(\nabla_S)|_{N_{T,S}}= (\pi_{T,S})^*\nabla_T  \quad\text{whenever
$T\prec S$}. \qquad\qquad
$$
\enddemo

 One sees that (4.1.1) and (4.3.1) imply that a $\Cal
C$-controlled connection on $E$ defines a usual connection on
$E|_S$ for every $S\in\Cal S$.  The next observation is evident
from the definition:

\proclaim
{(4.3.2) Lemma} The curvature form $\Theta\in j_*(A^2_{Y\!^\circ}\otimes\roman{
End}(E|_{Y\!^\circ}))$
of a $\Cal C$-controlled
connection $\nabla$ lies in $\Cal A^2_{Y,\Cal C}\otimes\roman{End}(E)$. \sq
\endproclaim

It is also obvious that $\Cal A^\bullet_{Y,\Cal C}$ is closed under
exterior
multiplication.  One can thus define for each $k$ the Chern form
$c_k(E,\nabla)$,
a closed {\it $\Cal C$-controlled} $2k$-form on $Y$, by the usual formula:
$$
c_k(E,\nabla) = P_k(\Theta,\dots,\Theta),
$$
where $P_k$ is the appropriate invariant polynomial of degree $k$.  By
(4.1.8),
$c_k(E,\nabla)$ defines a cohomology class in $H^{2k}(Y)$.
\proclaim
{(4.3.3) Proposition} i) Every $\Cal C$-controlled vector bundle $E$ on $Y$
admits a
$\Cal C$-controlled connection.
\smallskip

ii) The cohomology class of $c_k(E,\nabla)$ in $H^{2k}(Y)$ is independent
of the
$\Cal C$-controlled connection $\nabla$ on $E$.
\endproclaim
\demo
{Proof}
Let $\{\phi_T,\phi_{T,S}:T\prec S\}$ be the data defining a $\Cal C$-controlled
vector bundle,
as in (4.2.2.1), and $N'_T\subset N_T$ a domain for the isomorphisms
involving
$E_T$.  For each $T$, let $\nabla^T$ be any smooth connection on $E_T$,
and
$(\pi_T)^*\nabla^T$ the pullback connection on $E|_{N'_T}$.  Then $\frak V =
\{N'_T:T\in\Cal S\}$
is an open covering of $Y$.  Apply (4.1.6) to get a $\Cal C$-controlled
partition of
unity $\{f_T\}$ subordinate to $\frak V$.  Then $\nabla = \sum_V f_V\nabla^V$
is a
$\Cal C$-controlled connection on $E$.  This proves (i).

The argument for proving (ii) is the standard one.  For two connections on a
smooth
manifold, such as $Y^\circ$, there is an identity:
$$
c_k (E, \nabla_1) - c_k (E, \nabla_0) = d\eta_k,\tag 4.3.4
$$
where
$$
\eta_k = k \int^1_0 P_k(\omega, \Theta_t, \ldots, \Theta_t)dt,\tag 4.3.4.1
$$
$\omega = \nabla_1 - \nabla_0$, $\nabla_t = (1-t)\nabla_0 + t\nabla_1$, and
$\Theta_t$
denotes the curvature of $\nabla_t$.
Now, if
$\nabla_0$ and $\nabla_1$ are both $\Cal C$-controlled, one sees easily that
$\omega$ and
$\nabla_t$ are likewise, and then so is $\eta_k$. It follows that (4.3.4) is
an identity in
$\Cal A^\bullet_{Y,\Cal C}$, giving (ii).
\enddemo

We have been leading up to the following:
\proclaim
{(4.3.5) Theorem} Let $E$ be a $\Cal C$-controlled vector bundle on the
stratified space
$Y$. Then the
cohomology
class in {\rm (4.3.3,\,ii)} gives the topological Chern class of $E$ in
$H^{2k}(Y)$;
in particular, it is independent of the choice of $\Cal C$.
\endproclaim
\demo
{Proof} This argument, too, follows standard lines.  We start by proving the
assertion
when $E$ is a line bundle $L$.  On $Y$, there is the short exact exponential
sequence (of sheaves):
$$
0\to \ZZ_Y\to \Cal A^0_{Y,\Cal C}\to (\Cal A^0_{Y,\Cal C})^*\to 1.\tag 4.3.5.1
$$
The Chern class of $L$, $c_1 (L)$, is then the image of any controlled Cech
cocycle
that determines $L$, under the connecting homomorphism
$$
H^1(Y,(\Cal A^0_Y)^*)\longrightarrow H^2(Y,\ZZ).\tag 4.3.5.2
$$

To prove the theorem for line bundles, it is convenient to work in the double
complex $\frak
C^\bullet (\Cal A^\bullet_{Y,\Cal C})$,
where $\frak C^\bullet$ denotes Cech cochains.  It has differential $D= \delta
+ \sigma d$ (i.e., Cech
differential plus a sign $\sigma = (-1)^a$ times exterior derivative, where
$a$ is the Cech
degree).  On a sufficiently fine covering of $Y$ we have a cochain giving $L$,
$\xi\in
\frak C^1((\Cal A^0_{Y,\Cal C})^*)$ (if $e_\alpha$ is the specified
frame for
$L$ on the open subset $V_\alpha$ of $Y$, one has on $V_\alpha\cap V_\beta$
that
$e_\alpha = \xi_{\alpha \beta}e_\beta$), with $\delta\xi = 1$, the
connection
forms $\omega\in \frak C^0(\Cal A^1_{Y,\Cal C})$, and $\lambda = \log\xi$ in
$\frak
C^1(\Cal A^0_{Y,\Cal C})$.  We know by (4.3.5.2) above that $\delta\lambda$
gives
$c_1(L)$.  The change-of-frame formula for connections gives $\delta\omega +
d\lambda =0$.  Finally, the curvature
(for a line bundle)
is $\Theta = d\omega$, so we wish to show that $d\omega$ and $\delta\lambda$
are cohomologous
in the double complex.  By definition, $D\lambda = \delta\lambda - d\lambda$,
and
$D\omega = \delta\omega + d\omega = d\omega - d\lambda$.  This gives
$\delta\lambda
- d\omega = D(\lambda - \omega)$, and we are done.

To get at higher-rank bundles, we invoke a version of the splitting principle.
Let $p: \bold F(E)\to Y$
be the bundle of total flags for $E$.  As $E$ is locally the product of
$Y$ and a
vector space, $\bold F(E)$ is locally on $Y$ just $F_r\times Y$,
where $F_r$
is a (smooth compact) flag manifold.  As such, $\bold F(E)$ is a stratified
space that
is locally no
more complicated than $Y$;  we take as the set of strata $\widetilde{\Cal S} =
p^\-(\Cal S)=
\{p^\- S: S\in \Cal S\}$. For weak control data, we deduce it from $\Cal C$
in the same
way it is done for $E$ (see (4.2.2.2)): we take $N_{\bold F_T} = \bold F(E|_{
N'_T})$, and
use the natural projection $\bold F(E|_{N'_T})\to \bold F(E_T)$ induced by
(4.2.2.1).

It is standard that the vector bundle $p^*\!E$ on $\bold F(E)$ decomposes
(non-canonically) into a
direct sum of line bundles: $p^*\!E = \bigoplus_{1\leq j\leq r}\Lambda_j$.
($p^*\!E$
is canonically filtered: $\Lambda_1=F_1\subset F_2\dots F_r=p^*\!E$, with
$\Lambda_j\simeq F_j/F_{j-1}$.)
To obtain this, one starts by taking $\Lambda_1$ to be the line bundle given
at each
point of $\bold F(E)$ by the one-dimensional subspace from the corresponding
flag.  Then,
one splits the exact sequence
$$
0\to\Lambda_1\to p^*\!E\to p^*\!E/\Lambda_1\to 0,\tag 4.3.5.3
$$
using a {\it controlled metric} on $E$.  By that, we mean a metric that is a
pullback
via the isomorphisms (4.2.2.1); these can be constructed by the usual patching
argument,
using controlled partitions of unity (4.1.6). One obtains $\Lambda_j$, for
$j>1$, by
recursion.  We need a little more than that:
\proclaim{(4.3.6) Lemma}\,\,
i) The vector bundle $p^*\!E$ is, in a tautological way, a controlled vector
bundle on $\bold F(E)$.
\smallskip

ii) The line bundles $\Lambda_j$ are controlled subbundles of $p^*\!E$.
\endproclaim
\demo{Proof} We have that $p^*\!E = E\times_Y\bold F$, and its strata are
$(p^*\!E)_{{\bold F}_T} = E_T\times_T\bold F_T$, for all $T\in\Cal S$.  There
is natural weak control data for $p^*\!E$ that we now specify.
By construction, we have a retraction
$$
\pi_{(p^*\!E_T)}:(p^*\!E)|_{N_{\bold F_T}}=(p^*\!E)|_{\bold F|_{N_T}}\simeq
E|_{N_T}
\times_{N_T}\bold F|_{N_T}\to E_T\times_T\bold F_T=(p^*\!E)_{{\bold F}_T},
\tag 4.3.6.1
$$
induced by the weak control data for $E$ (and thus also $\bold F$), and
likewise for
the restriction to $(p^*\!E)|_{N_{\bold F_T,\,\bold F_S}}$, when $S\succ T$.
These provide $\phi_{\bold F_T}$ and $\phi_{\bold F_T,\,\bold F_S}$ (from
(4.2.2.1)) respectively for $p^*\!E$, and (i) is proved.

We show that $\Lambda_1$ is preserved by $\phi_{\bold F_T}$ and $\phi_{\bold
F_T,\,\bold F_S}$.  (As before, we explain this only for the former, the other
being its restriction to the strata.)
Let $p_T:\bold F_T\to T$ denote the restriction of $p$ to $\bold F_T$.
Since $p_T$ gives the flag manifold bundle associated to $E_T$, $(p_T^*\!E_T)$
contains a tautological line bundle, which we call $\Lambda_{1,T}$.
We have
from the control data that $(\Lambda_1)|_{N_{\bold F_T}}\simeq (\phi_{\bold
F_T})^*\!\Lambda_{1,T}$.

We claim further that (4.3.6.1) takes $\Lambda_1|_{N_{\bold F_T}}$ to
$\Lambda_{1,T}$,
as desired.
The explicit formula for (4.3.6.1), obtained by unwinding the fiber products,
is as follows.  Let $e$ be in the vector
space
$E_{T,t}$, the fiber of $E_T$ over $t\in T$, and $f$ a point the flag
manifold
of $E_{T,t}$.  Also, let $n\in (\pi_T)^\-(t)$.  Then $\pi_{(p^*\!E_T)}(e,f,n)
=(e,f)$,
which implies (ii) for $j=1$.  The assertion for $j>1$ is obtained recursively.
\sq
\enddemo

We return to the proof of (4.3.5).  Let $\nabla_0$ be the direct sum of
$\widetilde{
\Cal C}$-controlled connections on each $\Lambda_j$; and take
$\nabla_1 = p^*\nabla$, where $\nabla$ is a $\Cal C$-controlled connection on
$E$.
Both $\nabla_0$ and $\nabla_1$ are $\widetilde{\Cal C}$-controlled connections
on $\bold F(E)$.
By construction, $c_k(p^*\!E,\nabla_0)$ represents the $k$-th Chern class of
$p^*\!E$
in $H^{2k}(\bold F(E) )$.  We then apply (4.3.3,\,ii) to obtain that
$c_k(p^*\!E,
\nabla_1) = p^*c_k(E,\nabla)$ represents $p^*c_k(E)\in H^{2k}(\bold F(E) )$.
Since
$p^*: H^{2k}(Y )\longrightarrow H^{2k}(\bold F(E) )$ is injective, it
follows that
$c_k(E,\nabla)$ represents $c_k(E)$ in $H^{2k}(Y )$, and (4.3.5) is proved.
\sq
\enddemo
\bigskip

\centerline{\bf 5. Proofs of Theorem 2 and Conjecture B}
\medskip

In this section, we apply the methods of \S 4 in the case $Y=M^\RBS_\Gamma$.
\medskip

{\bf (5.1)} {\it Control data for a manifold-with-corners}.  Let $Y$ be a
manifold-with-corners,
with its open faces as strata.  For each codimension one boundary stratum $S$,
let
$$
\overline\phi_S: [0,1]\times\overline S\to Y
$$
define the collar $\overline N_S$ of $\overline S$
in $Y$,
so that $\{0\}\times \overline S$ is mapped identically onto $\overline S$.
This determines
partial control data (that is, without distance functions) for $Y$ as follows.

As $N_S$, one takes $\overline\phi([0,1)\times S)$, and as $\pi_S$ projection
onto $S$.
For a general boundary stratum $T$, write
$$
\overline T =\bigcap\{\overline S: S\text{ of codimension one, }T\prec S\}.
$$
Let $\overline N_T = \bigcap\{\overline N_S: S\text{ of codimension one, }
T\prec S\}$;
given the $\overline\phi_S$'s above, this set is canonically diffeomorphic
to $[0,1]^r
\times T$, where $r$ is the codimension of $T$.  Then $N_T$ is the subset of
$\overline
N_T$ corresponding to $[0,1)^r\times T$, in which terms $\pi_T$ is simply
projection
onto $T$.
\medskip

{\bf (5.2)} {\it Compatible control data.} The existence of natural (partial)
control data
for $M^\RBS_\Gamma$ is, in essence, well-known, as is compatible control data
for $M^\BB_\Gamma$
in the Hermitian case.  We give a brief presentation of that
here.
This will enable us to determine that Conjecture B is true.

The relevant notions are variants of (4.1.2).
\proclaim{(5.2.1) Definition {\rm (see [GM2,\,1.6])}} Let $Y$ and $Y'$ be
stratified
spaces.  A proper smooth mapping
$f:Y\to Y'$ is said to be {\rm stratified} when the following two conditions
are satisfied:
\smallskip

i) If $S'$ is a stratum of $Y'$, then $f^\-(S')$ is a union of connected
components of
strata of $Y$;
\smallskip

ii) Let $T\subset Y$ be a stratum component as in (i) above. Then
$f|_T:
T\to S'$ is a submersion.
\endproclaim

\noindent It follows that a stratified mapping $f$ is, in particular, open.
We assume
henceforth, and without loss of generality, that all strata are connected.

\proclaim{(5.2.2) Definition {\rm (cf.~[V1:\,1.4])}} Let $f:Y\to Y'$ be a
stratified
mapping, with (weak) control data $\Cal C$ for $Y$, and $\Cal C'$ for $Y'$.
We say
that $f$
is {\rm weakly controlled} if for each stratum $S$ of $Y$, the equation
$\pi_{S'}\circ f = f\circ\pi_{S}$
holds in some neighborhood of $S$ (here $f$ maps $S$ to $S'$).
\endproclaim

\noindent (5.2.3) {\it Remark.} Note that there is no mention of distance
functions in (5.2.2).
This is intentional, and is consistent with our stance in (4.1).
\smallskip

\proclaim{(5.2.4) Lemma} Let $f:Y\to Y'$ be a stratified mapping.  Given
partial control
data $\Cal C$ for $Y$, there is at most one system of germs of partial control
data $\Cal C'$
for $Y'$ such that $f$ becomes weakly controlled. Such $\Cal C'$ exists if and
only if
for all strata $S$ of $Y$, there is a neighborhood of $S$ (contained in $N_S$)
in
which $f(y)=f(z)$ implies $f(\pi_S(y))=f(\pi_S(z))$. \sq
\endproclaim

When the condition in (5.2.4) is satisfied, one uses the formula
$$
\pi_{S'}(f(y)) = f(\pi_{S}(y))
$$
to define $\Cal C'$, and we then write $\Cal C' = f_*\Cal C$.  In the usual
manner, the mapping
$f$ determines an equivalence relation on $Y$, viz., $y\sim z$ if and only if
$f(y)=
f(z)$.  The condition on $\Cal C$ thereby becomes:
$$
y\sim z\quad \Rightarrow\quad \pi_S(y)\sim \pi_S(z)\qquad\text{(near $S$)}.
\tag 5.2.5
$$

We will use the preceding for the stratified mappings $M^{BS}_\Gamma\to
M_\Gamma^\RBS$
in general, and $M_\Gamma^\RBS\to M_\Gamma^\BB$ in the Hermitian case.
The reason
for bringing in $M^{BS}_\Gamma$ is that it is a manifold-with-corners,
and it also has natural partial control data.

The boundary strata of $\Xbar$, the universal cover of $M_\Gamma^{BS}$, are
the sets $e(P)$
of (1.3), as $P$ ranges over all rational parabolic subgroups of $G$.  Those
of
$M_\Gamma^{BS}$
itself are the arithmetic quotients $e'(P)$ of $e(P)$, with $P$ ranging over
the finite set
of $\Gamma$-{\it conjugacy classes} of such $P$.  There are projections $X\to
X/A_P = e(P)$,
defined by collapsing the orbits of the geodesic action of $A_P$ to points.
This extends
to a $ P$-equivariant smooth retraction (geodesic projection), given in
(1.3.3):
$$
X(P) @>\pi_P >> X(P)/A_P = e(P).\tag 5.2.6
$$

Recall from (1.8) that there is a neighborhood of $\overline{e'(P)}$
in $M_\Gamma^{BS}$
on which geodesic projection onto $\overline{e'(P)}$,
induced by
(5.2.6), is defined.  We take the restriction of this geodesic projection
over $e'(P)$
as the definition of $\pi_P$ in our partial control data $\Cal C$ for
$M_\Gamma^{BS}$.
We have been leading up to:
\proclaim{(5.2.7) Proposition}
The quotient mapping $M^{BS}_\Gamma\to M_\Gamma^\RBS$ satisfies {\rm (5.2.5)}.
\endproclaim
\demo{Proof}
The mappings $e(P)\to e(P)/U_P$, as $P$ varies, induce the mapping $M^{
BS}_\Gamma
\to M_\Gamma^\RBS$.  It is a basic fact ([BS,\,4.3]) that for $Q\subset P$,
$A_Q\supset
A_P$ and $\pi_Q\circ\pi_P=\pi_Q$.  This gives $\overline{e(P)}\to
\overline{e(P)}/U_P$.  Since it is also the case that $Q
\subset P$
implies $U_Q\supset U_P$, we see that (5.2.5) is satisfied. \sq
\enddemo

Only a little more complicated is:
\proclaim{(5.2.8) Proposition}
In the Hermitian case, the quotient mapping $M_\Gamma^\RBS\to M_\Gamma^\BB$
satisfies {\rm (5.2.5)}.
\endproclaim
\demo{Proof}
When the symmetric space $X$ is Hermitian, the $P$-stratum of $X^\RBS$,
for each $P$,
decomposes as a product:
$$
e(P)/U_P\simeq X_{\ell,P}\times X_{h,P}.\tag 5.2.8.1
$$
This is induced by a decomposition of reductive algebraic groups over $\QQ$:
$$
P/U_P = G_{\ell,P}\cdot G_{h,P}
$$
(cf.~[Mu,\,p.\,\,254]).   Fixing $G_h$, one sees that the set of $Q$ with
$G_{h,Q}= G_h$
(if non-empty) is a lattice, whose greatest element is a maximal parabolic
subgroup $P$ of $G$.  The lattice is then canonically
isomorphic
to the lattice
of parabolic subgroups $R$ of $G_{\ell,P}$, whereby $Q/U_Q\simeq (R/U_R)\times
G_{h,P}$.
(Thus $G_{h,Q}= G_{h,P}$.  In the language of [HZ1,(2.2)] such $Q$ are said
to be {\it
subordinate to} $P$.)

The mapping $M_\Gamma^\RBS\to M_\Gamma^\BB$ is induced, in terms of (5.2.8.1),
by
$$
e(Q)/U_Q\to X_{h,Q},\tag 5.2.8.2
$$
for all $Q$;
perhaps more to the point, the terms can be grouped by lattice, yielding
$$
\overline X_{\ell,P}\times X_{h,P}\rightarrow X_{h,P}\hookrightarrow
(X_{h,P})^{BB}\tag 5.2.8.3
$$
for $P$ maximal (see [GT,\,2.6.3]).  One sees that (5.2.5) is satisfied. \sq
\enddemo
We have thereby reached the conclusion:
\proclaim{(5.2.9) Corollary}  The natural partial control data for
$M_\Gamma^{BS}$
induces compatible partial control data for $M_\Gamma^\RBS$ and $M_\Gamma^
\BB$. \sq
\endproclaim
\medskip

{\bf (5.3)} {\it Conjecture B.} Let $\CE_\Gamma$ be a homogeneous vector
bundle on
$M_\Gamma$, and $\CE_\Gamma^\RBS$
its extension to $M_\Gamma^\RBS$ from [GT] that was reconstructed in our
\S 1. We select
as partial control data $\Cal C$ for $M_\Gamma^\RBS$ that given in (5.2.9).
It is essential
that the following hold:
\proclaim{(5.3.1) Proposition} $\CE_\Gamma^\RBS$ is a controlled vector bundle
on $M_\Gamma^\RBS$,
with the $\widetilde\pi_P$'s of {\rm (1.3.3)} providing the weak control data.
\endproclaim
\demo{Proof}  This is
almost
immediate from the construction in \S 1.  Recall that the weak control data
for $M_\Gamma^\RBS$
consists of the geodesic projections $\pi_P$, defined in a neighborhood of
$\widehat M_P$.
The vector bundle $\CE_\Gamma^\RBS$ also gets local geodesic projections
$\widetilde\pi_P$,
induced from those of $\CE^\BS$, that are compatible with those of $M_\Gamma^
\RBS$ because
of (1.2).  The same holds within the strata of these spaces, by (1.6).  We
see that
the criterion of (4.2.3) is satisfied. \sq
\enddemo
We proceed with a treatment of $\nabla^{\text{GP}}$, the connection on
$\CE_\Gamma$ constructed
in [GP].  For each maximal $\QQ$-parabolic subgroup of $G$, let $M_P$ be the
corresponding
stratum of $M^\BB$; it is a locally symmetric variety for the group $G_{h,P}$.
We also use
``$P$'' to label the strata: thus, we have for (4.1.2), $\pi_P:N_P\to M_P$,
etc.
Then $\nabla^{\text{GP}}$
can be defined recursively, starting from the strata of lowest dimension
($\QQ$-rank zero),
and then increasing the $\QQ$-rank by one at each step.

There is, first, the equivariant {\it Nomizu connection} for homogeneous
vector
bundles, whose definition we recall.  Homogeneous vector bundles are
associated
bundles of the principal $K$-bundle:
$$
\kappa: \Gamma\backslash G\longrightarrow M_\Gamma.\tag 5.3.2
$$
When we write the Cartan decomposition $\g = \k \oplus\p$, we note that
(5.3.2)
has a natural
equivariant connection whose connection form lies in the vector space
$\text{Hom} (\g,\k)$;
it is given by the projection of $\g$ onto $\k$ (with kernel $\p$).  This
is known as
the Nomizu connection.  The homogeneous vector bundle $\Cal E_\Gamma$
on $M_\Gamma$
is associated to the principal bundle (5.3.2) via the representation
$K\to \text{GL}(E)$.
The connection induced on $\Cal E_\Gamma$ via $\k\to\frak{gl}(E)=\text{End}
(E)$ is also
called the
Nomizu connection (of $\Cal E_\Gamma$), and will be denoted $\nabla^{N\!o}$;
its connection form is denoted
$\theta\in\g^*\otimes \text{End}(E)$.
A $K$-frame for
$\Cal E_\Gamma$ on an open subset $O\subset M_\Gamma$ is given by a smooth
cross-section
$\sigma: O\to \kappa^\- (O)$ of $\kappa$; the resulting connection matrix is
the pullback
of $\theta$ via $\sigma^*$, an element of $A^1(O,\text{End}(E))$.

With that stated, we can start to describe $\nabla^{\text{GP}}$.  For any
maximal
$\QQ$-parabolic
$P$, one will be taking expressions of the form
$$
\nabla^P = \psi^P_P\nabla^{P,N\!o} + \sum_{Q\prec P}\psi^Q_P\Phi^*_{Q,P}
(\nabla^Q)
\tag 5.3.3
$$
[GP,\,11.2]. Here, $\nabla^{P,N\!o}$ is the Nomizu connection for the
homogeneous vector
bundle
on $M_P$ determined by the restriction $K_{h,P}\hookrightarrow K\to\text{GL}(E)$,
and the functions
$\{\psi^Q_P: Q\preceq P\}$ form a partition of unity on $M_P$ of a selected
type,
given in [GP,\,3.5,\,11.1.1]; the function $\psi^Q_P$ is a cut-off function
for a large
relatively compact open subset $V_{Q,P}$ of $M_Q$ in $M^*_P$, with
$$
\bigcup_{Q\preceq P}V_{Q,P} = M^*_P
,
$$
and can be taken to be
supported
inside the neighborhood $N_{Q,P}$ of the partial control data when $Q\ne P$.

Next, $\Phi^*_{Q,P}$ indicates the process of {\it parabolic induction}
from $M_Q$
to $N_{Q,P}$, by means of $\pi_{Q,P}$.  It is defined as follows.  Fix a
maximal parabolic
$Q$ and a
representation
$K\to \text{GL}(E)$.  The latter restricts, of course, to $K_Q=K_{h,Q}\times
K_{\ell,Q}$,
but through the Cayley transform, this actually extends to a representation
$\lambda$ of $K_{h,Q}\times
G_{\ell,Q}$.
That allows one to define an action of {\it all of} $Q$ on $\Cal E_{h,Q}$
[GP,\,10.1],
which induces a $Q$-equivariant mapping
$$
\Cal E = Q\times_{K_Q}E\longrightarrow G_{h,Q}\times_{K_{h,Q}}E = \Cal
E_{h,Q},\tag 5.3.4
$$
given by $(q,e)\mapsto (g_h,\lambda(g_\ell)e)$ for $q=ug_hg_\ell\in U_Q
G_{h,Q}G_{\ell,Q}=Q$.
That in turn defines a $U_Q$-invariant isomorphism of vector bundles
homogeneous
under $Q$:
$$
\Cal E\simeq \pi_Q^*(\Cal E_{h,Q}).
$$

One then takes $\nabla^{\text{GP}}$ to be $\nabla^G$ in (5.3.3).
Given any connection
$\nabla^{h,Q}$
on $\Cal E_{h,Q}$, the pullback connection $\nabla = \pi_Q^*(\nabla^{h,Q})$
satisfies
the same relation for its curvature form, viz.,
$$
\Theta(\nabla) = \pi_Q^*\Theta(\nabla^{h,Q}).\tag 5.3.5
$$
It follows that the Chern forms of $\nabla^{\text{GP}}$ are controlled on
$M^\BB_\Gamma$
[GP,\,11.6].  On the other hand, the connection itself is not.  To proceed,
weaker
information about $\nabla^{\text{GP}}$ suffices:
\proclaim{(5.3.6) Proposition} With an appropriate choice of the functions
$\psi^Q_P$,
the connection $\nabla^{\text{GP}}$ {\rm is} a controlled connection
when viewed on
$M^\RBS_\Gamma$.
\endproclaim
\demo{Proof} TRBS1-2
his is not difficult.  Recall from (4.3.1) that the issue is
the existence
of local frames at each point of $M^\RBS_\Gamma$, with respect to which the
connection
matrix is controlled.  For each rational parabolic subgroup $Q$ of $G$, we
work in the
corner $X(Q)$.  By (5.3.4), one gets local frames for $\Cal E(Q)$, the
restriction
of $\CEbar$ to $X(Q)\subset\Xbar$, from local frames for $ \Cal E_{h,Q}$.  We
can write $\Cal E(Q)$ as:
$$
\Cal E(Q)\simeq U_Q\times \Abar_Q\times M_Q\times_{K_Q} E.\tag 5.3.6.1
$$
This also provides good variables for calculations.  We note that
$\Phi^*_{Q,P}$
is independent of the $U_Q$-variable.  Likewise, $\psi^Q_P$ can be chosen
to be a
function of only $(a,mK_Q)$, constant on the compact nilmanifold fibers $N_P$
(i.e., the
image of the $U_P$-orbits).  It follows by induction that $\nabla^{\text{GP}}$
is
controlled on $M^\RBS_\Gamma$. \sq
\enddemo

As we said, the Chern forms of $\nabla^{\text{GP}}$
are
controlled differential forms for $M^\BB_\Gamma$, so are {\it a fortiori}
controlled
for $M^\RBS_\Gamma$.  It follows from (4.3.5) that
\proclaim{(5.3.7) Proposition}
$c_\b (\CE_\Gamma
^\RBS,\nabla^{\text{GP}})$ represents $c_\b (\CE_\Gamma^\RBS)\in H^\b
(M^\RBS_\Gamma )$. \sq
\endproclaim
Thereby, Conjecture B is proved.
\bigskip

{\bf (5.4)} {\it Theorem 2.}
Let $\nabla^\control$ be any $\Cal C$-controlled connection on $\CE^
\RBS_\Gamma$,
and $\nabla^{N\!o}$ the equivariant Nomizu connection on $\CE_\Gamma$. By
(4.3.5) we
know that
$c_k(\nabla^\control)$ represents $c_k(\CE_\Gamma^\RBS)$; we want to conclude
the same
for $c_k(\nabla^{N\!o})$.  Toward that,
we recall the standard identity on $M$ satisfied by the Chern forms:
$$
c_k(\nabla^{N\!o})-c_k(\nabla^\control) = d\eta_k,\tag 5.4.1
$$
which is a case of (4.3.4).  The following is straightforward:

\proclaim{(5.4.2) Lemma} i) $\Cal A^\bullet_{M^\RBS_\Gamma\!,\,\Cal C}$ is
contained
in $\Cal A_{(\infty)}^\bullet (M^\RBS_\Gamma)$.

ii) A $G$-invariant form on $M_\Gamma$ is $L^\infty$.
\endproclaim
\demo{Proof} In terms of (3.1.1), a controlled differential form on $M^\RBS_
\Gamma$ is
one that is, for each given $P$, pulled back from $V\subset \widehat M_P$.
Such forms
are trivially weighted by $A_P^+$ in the
metric (3.1.2).
It follows that a controlled form is locally $L^\infty$ on $M^\RBS_\Gamma$.
This proves (i).
As for (ii), an invariant form has constant length, so is in particular
$L^\infty$. \sq
\enddemo

\proclaim{(5.4.3) Proposition}The closed forms $c_k(\nabla^{N\!o})$ and
$c_k(\nabla^
\control)$ represent the same class in $H^{2k}_{(\infty)}(M_\Gamma)$.
\endproclaim
\demo{Proof}
Since $M^\RBS_\Gamma$ is compact, a global controlled form
on $M^\RBS_
\Gamma$ is globally $L^\infty$.  As such, (5.4.2) gives that the Chern forms
for
both $\nabla^
\control$ and $\nabla^{N\!o}$ are in the complex $L^\bullet_{(\infty)}(M_
\Gamma)$.
It remains
to verify that $\eta_k$ in (5.4.1) is likewise $L^\infty$, for then the
relation
(5.4.1) holds in the $L^\infty$ de Rham complex $A^\bullet_{(\infty)}(M_
\Gamma)$,
so $c_k(\nabla^{N\!o})$ and $c_k(\nabla^\control)$ are cohomologous in
the $L^\infty$ complex.

By (4.3.4.1), it suffices to check that the difference $\omega = \nabla^{
N\!o} - \nabla^
\control$ is $L^\infty$.  That can be accomplished by taking the difference
of connection
matrices with respect to the same local frame of $\Cal E^\RBS$,
and for
that purpose we use, for each $Q$, frames pulled back from $\widehat M_Q$.
For that,
it is enough to verify the boundedness for $\omega$ in a neighborhood
of every point
of the boundary of $M^\RBS_\Gamma$, and we may as well calculate on
$X^\RBS$.

Consider a point in the $Q$-stratum $X_Q$ of $X^\RBS$.  As in  (3.1.1), we can
take as
neighborhood
base, intersected with $X$, sets that decompose with respect to $Q$ as
$$
N_Q\times A_Q^+ \times V,\tag 5.4.3.1
$$
with $V$ open in $X_Q$. In these terms, $\pi_Q$ is just projection onto $V$.
As in
(3.1.2) and (3.1.4), we use as coordinates $(u_\alpha,a,v)$.  We also
decompose
(see the end of (1.1)),
$$\CE\simeq Q\times_{K_Q} E\simeq U_Q\times A_Q^+\times X_Q\times_{
K_Q} E.\tag 5.4.3.2
$$
We obtain a canonical isomorphism $\CE\simeq\pi_Q^*\CE_Q$, with $\CE_Q$
a homogeneous
vector bundle on $X_Q$.  By (5.4.2,\,ii), the connection matrix of a connection
that
is pulled back from $X_Q$, with respect
to a
local frame
pulled back from $X_Q$ is $L^\infty$, so we wish to do the same for the
Nomizu connection.

First, we have:
\proclaim{(5.4.3.3) Lemma} Let $\widehat Q = Q/A_QU_Q$ and consider the
diagram
$$\CD
Q & @>>> & \widehat Q\\
@VVV & &  @VVV\\
X & @>\pi_Q >> & X_Q.\endCD
$$
Then $Q\simeq \widehat Q\times_{X_Q} X$, the pullback of $\widehat Q$ with
respect to
$\pi_Q$.
\endproclaim
\demo{Proof}
We note that both $\widehat Q$ and $Q$ are exhibited as principal
$K_Q$-bundles.  To
prove our assertion,
it is simplest use the Langlands decomposition (of manifolds) $Q\simeq\widehat
Q\times A_Q \times U_Q$
to yield the decomposition $X\simeq X_Q\times A_Q \times U_Q$ (cf.~(5.4.3.1)).
Then
$$
\widehat Q\times_{X_Q} X\simeq \widehat Q\times A_Q \times U_Q\simeq Q.
\quad\square
$$
\enddemo

It follows that if $\sigma: O\subseteq X_Q\to \widehat Q$ gives a local
$K_Q$-frame, then
$\widetilde
\sigma : \pi_Q^\- (O)\subseteq X\to Q\simeq \widehat Q\times_{X_Q} X$, defined
by $\widetilde
\sigma(x)= (\sigma(\pi_Q(x)),x)$, gives the pullback frame $\pi_Q^*\sigma$.
In other
words, $\pi_Q^*\sigma$ takes values in the principal $K_Q$-bundle $Q\to X$
that is the
restriction of structure group of (5.3.2) from $K$ to $K_Q$.

Let $\nabla^{N\!o}$ be the Nomizu connection on $\CE$.  Recall that this is
determined by
$$
T_X @> \widetilde\sigma_* >> \g @>>> \k @>>> \roman {End}(E).\tag 5.4.3.4
$$
As such, $\nabla^{N\!o}$ is not a $K_Q$-connection.  However,
a frame
for the restriction to $X$ of the canonical extension $\CEbar$ can be taken to
be of the
form $\widetilde\sigma$ as above (cf.~(1.10)).  It follows that for
$x\in\pi_Q^\- (O)$,
the Nomizu connection is given by
$$
T_{X,x} \overset \widetilde\sigma_*\to\hookrightarrow \q @>>> \k @>>>
\roman {End}(E),\tag 5.4.3.5
$$
where $\q$ denotes the Lie algebra of $Q$.  This is a mapping that is of
constant
norm along the fibers of $\pi_Q$.  It follows that the connection matrix is
$L^\infty$.
Therefore, we have:

\proclaim{(5.4.3.6) Proposition}
The connection difference $\omega$ is $L^\infty$. \sq
\endproclaim
\noindent This finishes the proof of (5.4.3).

\demo{{\rm (5.4.4)} Remark} The reader may find it instructive to compare,
in the case of
$G=SL(2)$, the
above argument
to the one used in [Mu,\,pp.~259--260].  The two discussions, seemingly quite
different, are
effectively the same.
\enddemo

We now finish the proof of Theorem 2 by demonstrating:

\proclaim{(5.4.5) Proposition} $c_k(\nabla^{N\!o})$ and $c_k(\nabla^
\control)$
represent the same class in $H^{2k}(M^\RBS_\Gamma)$.
\endproclaim
\demo{Proof} Because $M_\Gamma$ has finite volume, there is a canonical
mapping
$$
H^\b _{(\infty)}(M_\Gamma)\to H^\b _{(p)}(M_\Gamma)
$$
for all $p$ (see (3.3.1)).  It follows from (5.4.3) that $c_k(\nabla^{N\!o})$
and $c_k(\nabla^
\control)$ represent
the same class in $H^{2k}_{(p)}(M_\Gamma)$ for all $p$.
Taking
$p$ sufficiently large, we apply Theorem 1 (i.e., (3.1.11)) to see that
$c_k(\nabla^
{N\!o})$ and
$c_k(\nabla^\control)$ represent the same class in $H^{2k} (M^\RBS_\Gamma )$.
\sq
\enddemo
\bigskip\medskip

\centerline{\bf References}
\medskip

\noindent [BB] Baily, W., Borel, A., {\it Compactification of arithmetic
quotients
of bounded symmetric domains}.  Ann.~of Math.~{\bf 84} (1966), 442--528.

\smallskip
\noindent [BS] Borel, A., Serre, J.-P., {\it Corners and arithmetic
groups}.
Comm.~Math.~Helv.~{\bf 4} (1973), 436--491.
\smallskip

\noindent [GHM] Goresky, M., Harder, G., MacPherson, R., {\it Weighted
cohomology}.
Invent.~Math.~{\bf 116} (1994), 139--213.
\smallskip

\noindent [GM1] Goresky, M., MacPherson, R., {\it Intersection homology, II}.
Invent.~Math.~{\bf 72} (1983), 77--129.
\smallskip

\noindent [GM2] Goresky, M., MacPherson, R., Stratified Morse Theory.
Springer-Verlag, 1988.
\smallskip

\noindent [GP] Goresky, M., Pardon, W., {\it Chern classes of modular
varieties}, 1998.
\smallskip

\noindent [GT] Goresky, M., Tai, Y.-S., {\it Toroidal and reductive
Borel-Serre
compactifications of locally symmetric spaces}.
Amer.~J.~Math.~{\bf 121} (1999), 1095--1151.
\smallskip

\noindent [HZ1] Harris, M., Zucker, S., {\it Boundary cohomology of Shimura
varieties, II: Hodge theory at the boundary}.
Invent.~Math.~{\bf 116} (1994), 243--307.
\smallskip

\noindent [HZ2] Harris, M., Zucker, S., {\it Boundary cohomology of Shimura
varieties, III:
Coherent cohomology on higher-rank boundary strata and applications to
Hodge theory} (to appear, Mem.~Soc.~Math.~France).
\smallskip

\noindent [K] Kostant, B., Lie algebra cohomology and the generalized
Borel-Weil
theorem. Ann.~of Math.~{\bf 74} (1961), 329--387.
\smallskip

\noindent [L] Looijenga, E., {\it $L_2$-cohomology of locally symmetric
varieties.} Compositio Math.~{\bf 67} (1988), 3--20.

\noindent [M] MacPherson, R., {\it Chern classes for singular algebraic
varieties.}
Annals of Math.~{\bf 100} (1974), 423--432.
\smallskip

\noindent [Mu] Mumford, D.,  {\it Hirzebruch's proportionality theorem in the
non-compact
case}. Invent.~Math.~{\bf 42} (1977), 239--272.
\smallskip

\noindent [N] Nomizu, K., On the cohomology of compact homogeneous spaces of
nilpotent
Lie groups.  Ann.~of Math.~{\bf 59} (1954), 531--538.
\smallskip

\noindent [SS] Saper, L., Stern M., {\it $L_2$-cohomology of arithmetic
        varieties.}  Ann.~of Math.~{\bf 132} (1990), 1--69.
\smallskip

\noindent [V1] Verona, A., {\it Le th\'eor\`eme de de Rham pour les
pr\'estratifications
abstraites}.  C.~R.~Acad.~Sc.~Paris {\bf 273} (1971), 886--889.
\smallskip

\noindent [V2] Verona, A., {\it Homological properties of abstract
prestratifications}.
Rev.~Roum.~Math. Pures et Appl.~{\bf 17} (1972), 1109--1121.
\smallskip

\noindent [Z1]  Zucker, S.,  {\it $L_2$-cohomology of warped products and
arithmetic
groups.}  Invent.~Math.~{\bf 70} (1982), 169--218.
\smallskip

\noindent [Z2] Zucker, S., {\it $L_2$-cohomology and intersection homology
of locally
symmetric varieties, III.} In:
Hodge Theory: Proceedings Luminy, 1987, Ast\'erisque {\bf 179--180} (1989),
245--278.
\smallskip

\noindent [Z3] Zucker, S., {\it $L^p$-cohomology and Satake
compactifications}.
 In:
J.~Noguchi, T.~Ohsawa (eds.), {\it Prospects in Complex  Geometry: Proceedings,
Katata/Kyoto
1989}. Springer LNM {\bf 1468} (1991), 317--339.
\smallskip

\noindent [Z4] Zucker, S., {\it $L^p$-cohomology: Banach spaces and
homological
methods on Riemannian
manifolds}.
In: {\it Differential Geometry: Geometry in Mathematical Physics and Related
Topics},
Proc.~of Symposia in Pure Math.~{\bf 54} (1993), 637--655.
\smallskip

\noindent [Z5] Zucker, S., {\it On the boundary cohomology of locally
symmetric
varieties}.
Vietnam J.~ Math.~{\bf 25} (1997), 279--318, Springer-Verlag.

\end

Let $\nabla^{N\!o}$ be the Nomizu connection on $\CE$.  Recall that this is
determined by
$$
T_X @> \sigma_* >> \g\to\k\to \text{End}(E).\tag 5.4.3.4
$$
The above observation implies that this can be given instead by:
$$
T_X @> \sigma_* >> \q\to\widehat\q\to\k_Q\hookrightarrow\k\to \text{End}(E)
\tag 5.4.3.5
$$
({\it not} $\q\hookrightarrow\g\to\k$).  Here, $\q$ is the Lie algebra of $Q$,
$\widehat\q$ that
of $\widehat Q$, and $\k_Q$ that of $K_Q$.  It follows that the connection
matrix
of $\nabla^{N\!o}$ is $L^\infty$.  Therefore, we have: